\documentclass[leqno,11pt,a4paper]{amsart} 

\usepackage[top=28truemm,bottom=28truemm,left=38truemm,right=38truemm]{geometry}
\usepackage{fix-cm}
\usepackage[TU]{fontenc}
\makeatletter
\renewcommand\section{\@startsection{section}{1}{\z@}%
           {25\p@ \@plus 6\p@ \@minus 3\p@}%
           {10\p@ \@plus 6\p@ \@minus 3\p@}%
           {\fontsize{13pt}{0cm}\selectfont\bfseries\boldmath}}

\renewcommand\subsection{\@startsection{subsection}{2}{\z@}%
           {13\p@ \@plus 6\p@ \@minus 3\p@}%
           {6\p@ \@plus 6\p@ \@minus 3\p@}%
           {\fontsize{12pt}{0cm}}}
\renewcommand\subsubsection{\@startsection{subsubsection}{3}{\z@}%
           {12\p@ \@plus 6\p@ \@minus 3\p@}%
           {\p@}%
           {\normalfont\normalsize}}
\renewcommand{\paragraph}[1]{%
  \par
  \addvspace{\medskipamount}
  \noindent
  \textbf{#1\@addpunct{.}}\enspace\ignorespaces
}
\makeatother
\let\oldtocsection=\tocsection
\let\oldtocsubsection=\tocsubsection
\let\oldtocsubsubsection=\tocsubsubsection
\renewcommand{\tocsection}[2]{
\hspace{0em}\bfseries\oldtocsection{#1}{#2}}
\renewcommand{\tocsubsection}[2]{\hspace{1em}\small\oldtocsubsection{#1}{#2}}
\renewcommand{\tocsubsubsection}[2]{\hspace{2em}\small\oldtocsubsubsection{#1}{#2}}
\usepackage[x11names,dvipsnames]{xcolor}
\usepackage{graphicx}
\usepackage{url}
\usepackage[doi=false,isbn=false, sortcites=true,eprint=true, date=year, style=numeric,uniquename=init,giveninits=true,sorting=nty,url=true,maxnames=99,backref=true]{biblatex} 

\DeclareFieldFormat*{volume}{\mkbibbold{#1}} 
\DeclareFieldFormat*{title}{\mkbibquote{#1\adddot}}  
\DeclareFieldFormat*{booktitle}{\mkbibemph{#1}}  
\DeclareSourcemap{
  \maps[datatype=bibtex]{
    \map{
      \step[fieldsource=edition,
            match=\regexp{^(\d+)\.\s+(ed.|edition)$},
            replace=\regexp{$1}]
    }
  }
}
\DeclareFieldFormat*{journaltitle}{\mkbibemph{#1}} 
\DeclareFieldFormat[article,Article,book,inbook,incollection,inproceedings,patent,thesis,unpublished,misc]{year}{#1} %
\renewbibmacro{in:}{} %
\renewbibmacro*{issue+date}{%
  \printfield{issue}%
  \setunit*{\addspace}%
  \usebibmacro{date}%
  \newunit} %
\DeclareSourcemap{
	\maps[datatype=bibtex]{
		\map[overwrite=true]{
			\step[fieldsource=fjournal]
			\step[fieldset=journal, origfieldval]
		}
	}
}
\usepackage{xpatch}
\usepackage[hypertexnames=false]{hyperref}
\hypersetup{
    colorlinks=true,
    citecolor=blue,
    linkcolor=blue,
    urlcolor=teal,
    breaklinks=true
}
\usepackage{graphicx}
\usepackage{amsfonts}
\usepackage[capitalize,noabbrev,nameinlink]{cleveref}

\usepackage{amsthm}
\usepackage{thmtools}
\usepackage{thm-restate}
\usepackage{paralist}
\expandafter\let\expandafter\oldproof\csname\string\proof\endcsname
\let\oldendproof\endproof
\renewenvironment{proof}[1][\proofname]{%
  \oldproof[\bfseries\itshape #1] 
}{\oldendproof}


\theoremstyle{plain}
\newtheorem{theorem}{Theorem}[section]
\newtheorem{proposition}[theorem]{Proposition}
\newtheorem{lemma}[theorem]{Lemma}
\newtheorem{corollary}[theorem]{Corollary}

\theoremstyle{definition}
\newtheorem{definition}[theorem]{Definition}

\newtheorem{assumption}[theorem]{Assumption}

\theoremstyle{remark}
\newtheorem{remark}[theorem]{Remark}
\newtheorem{example}[theorem]{Example}

\usepackage{array}
\usepackage{latexsym}
\usepackage{subfiles}
\usepackage{optidef}
\allowdisplaybreaks[4]
\usepackage{stackengine}
\usepackage{empheq}

\usepackage{here}
\usepackage{longtable}
\usepackage{lipsum}
\usepackage{subcaption}
\usepackage{caption}
\usepackage{capt-of}
\usepackage{threeparttable}
\usepackage{upgreek}
\usepackage{boldline}
\usepackage{multirow}
\usepackage{here}
\usepackage[super]{nth}
\usepackage{svg}
\usepackage[inline]{enumitem}
\setlist[description]{%
  topsep=10pt,               
  itemsep=5pt,               
  font={\bfseries\rmfamily}, 
}
\newlist{myenum}{enumerate}{3}
\setlist[myenum,1]{label=\textbf{\arabic*.},
                   ref  =\textbf{\arabic*.}}
\setlist[myenum,2]{label=\textbf{(\alph*)},
                   ref  =\themyenumi\textbf{(\alph*)}}
\setlist[myenum,3]{label=\bfseries(\roman*),
                   ref  =\themyenumii\textbf{.(\roman*)}}

\usepackage{mathtools}
\usepackage{etoolbox}
\usepackage{breqn}
\usepackage{physics}
\usepackage{amssymb,amsmath}

\usepackage{lmodern}
\usepackage{mathrsfs}

\DeclareFontFamily{U}{matha}{\hyphenchar\font45}
\DeclareFontShape{U}{matha}{m}{n}{
<-6> matha5 <6-7> matha6 <7-8> matha7
<8-9> matha8 <9-10> matha9
<10-12> matha10 <12-> matha12
}{}
\DeclareSymbolFont{matha}{U}{matha}{m}{n}

\DeclareFontFamily{U}{mathx}{\hyphenchar\font45}
\DeclareFontShape{U}{mathx}{m}{n}{
<-6> mathx5 <6-7> mathx6 <7-8> mathx7
<8-9> mathx8 <9-10> mathx9
<10-12> mathx10 <12-> mathx12
}{}
\DeclareSymbolFont{mathx}{U}{mathx}{m}{n}

\DeclareMathDelimiter{\vvvert} {0}{matha}{"7E}{mathx}{"17}%

\DeclarePairedDelimiterX{\normiii}[1]
{\vvvert}
{\vvvert}
{\ifblank{#1}{\:\cdot\:}{#1}}
\usepackage{appendix}

\DeclareSymbolFont{EulerExtension}{U}{euex}{m}{n}
\DeclareMathSymbol{\euintop}{\mathop} {EulerExtension}{"52}
\DeclareMathSymbol{\euointop}{\mathop} {EulerExtension}{"48}

\numberwithin{equation}{section}

\allowdisplaybreaks[4]

\usepackage{makecell}
\usepackage{tikz}
\usetikzlibrary{calc}
\tikzset{
    ncbar angle/.initial=90,
    ncbar/.style={
        to path=(\tikztostart)
        -- ($(\tikztostart)!#1!\pgfkeysvalueof{/tikz/ncbar angle}:(\tikztotarget)$)
        -- ($(\tikztotarget)!($(\tikztostart)!#1!\pgfkeysvalueof{/tikz/ncbar angle}:(\tikztotarget)$)!\pgfkeysvalueof{/tikz/ncbar angle}:(\tikztostart)$)
        -- (\tikztotarget)
    },
    ncbar/.default=0.5cm,
}

\tikzset{square left brace/.style={ncbar=0.5cm}}
\tikzset{square right brace/.style={ncbar=-0.5cm}}

\tikzset{round left paren/.style={ncbar=0.5cm,out=120,in=-120}}
\tikzset{round right paren/.style={ncbar=0.5cm,out=60,in=-60}}
\usetikzlibrary{positioning, fit}   
\tikzset{block/.style={draw, thick, text width=2cm ,minimum height=1.3cm, align=center},   
line/.style={-latex}     
}

\crefname{theorem}{Theorem}{Theorems}
\crefname{definition}{Definition}{Definitions}
\crefname{problem}{Problem}{Problems}
\crefname{fact}{Fact}{Facts}
\crefname{proposition}{Proposition}{Propositions}
\crefname{lemma}{Lemma}{Lemmas}
\crefname{corolary}{Corolary}{Corolaries}
\crefname{conjecture}{Conjecture}{Conjectures}
\crefname{assumption}{Assumption}{Assumptions}
\crefname{enumi}{Assumption}{Assumptions}
\crefname{claim}{Claim}{Claims}

\crefname{remark}{Remark}{Remarks}
\crefname{example}{Example}{Examples}
\crefname{corollary}{Corollary}{Corollaries}
\crefname{subsubsection}{Subsubsection}{Subsubsections}
\crefname{subsection}{Subsection}{Subsections}
\crefname{section}{Section}{Sections}
\crefname{chapter}{Chapter}{Chapters}
\crefname{table}{Table}{Tables}
\crefname{figure}{Figure}{Figures}
\crefname{algorithm}{Algorithm}{Algorithms}

\crefname{myenumi}{item}{items}
\crefname{myenumii}{item}{items}
\crefname{myenumiii}{item}{items}

\renewcommand{\ref}{\cref}

\newcommand{\e}{\mathrm{e}}
\newcommand{\distance}{\mathsf{d}}

\newcommand{\Pcal}{\mathcal{P}}
\newcommand{\Ptwo}{\Pcal_2}
\newcommand{\Pac}{\Pcal_{\mathrm{ac},2}}

\newcommand{\Ccinf}{C_c^\infty}
\newcommand{\R}{\mathbb{R}}

\newcommand{\exR}{[-\infty,+\infty]}
\newcommand{\overR}{(-\infty,+\infty]}
\newcommand{\underR}{[-\infty,+\infty)}
\newcommand{\Z}{\mathbb{Z}}
\newcommand{\N}{\mathbb{N}}
\newcommand{\Y}{\mathcal{Y}}
\newcommand{\X}{\mathcal{X}}
\newcommand{\XY}{X \times Y}
\newcommand{\cyl}{\mathsf{Cyl}}

\newcommand{\Loneloc}{L^1_{\mathrm{loc}}}
\newcommand{\zeroinfty}{[0,+\infty)}

\newcommand{\relEnt}{\mathcal{H}}

\newcommand{\la}{\left\langle}
\newcommand{\ra}{\right\rangle}

\newcommand{\EVSI}[4]{\mathrm{EVIs}_{#4}(#1,#2,#3)}

\newcommand{\fsubd}{\underline{\partial}}
\newcommand{\fsupd}{\overline{\partial}}

\makeatletter
\newcommand{\pushright}[1]{\ifmeasuring@#1\else\omit\hfill$\displaystyle#1$\fi\ignorespaces}
\newcommand{\pushleft}[1]{\ifmeasuring@#1\else\omit$\displaystyle#1$\hfill\fi\ignorespaces}
\makeatother
\DeclareMathOperator{\cl}{\mathsf{cl}}

\DeclareMathOperator{\Domain}{\mathrm{Dom}}

\DeclareMathOperator*{\argmax}{\mathop{\rm arg~max}}
\DeclareMathOperator*{\argmin}{\mathop{\rm arg~min}}
\DeclareMathOperator*{\argminimax}{\mathop{\rm arg~minimax}}
\DeclareMathOperator*{\Lip}{\mathrm{Lip}}
\DeclareMathOperator*{\supp}{\mathrm{supp}}

\DeclareMathOperator{\Tan}{\mathrm{Tan}}

\def\Set#1{\Setdef#1\Setdef}
\def\Setdef#1|#2\Setdef{\left\{#1\,\;\mathstrut\vrule\,\;#2\right\}}%
\newcommand{\given}[2]{
  \left(#1 \;\middle\vert\; #2\right)
}
\renewcommand{\th}{%
    \ifmmode
        ^\mathrm{th}%
    \else%
        \textsuperscript{th}\xspace%
    \fi%
}
\makeatletter
\newcommand{\subalign}[1]{%
  \vcenter{%
    \Let@ \restore@math@cr \default@tag
    \baselineskip\fontdimen10 \scriptfont\tw@
    \advance\baselineskip\fontdimen12 \scriptfont\tw@
    \lineskip\thr@@\fontdimen8 \scriptfont\thr@@
    \lineskiplimit\lineskip
    \ialign{\hfil$\m@th\scriptstyle##$&$\m@th\scriptstyle{}##$\hfil\crcr
      #1\crcr
    }%
  }%
}
\makeatother

\newcommand{\Urd}{{\frac{ \dd}{\dd t}\!\!}^+}

\newcommand{\contiphi}{\varphi}
\newcommand{\phiX}{\psi_X}
\newcommand{\phiY}{\psi_Y}
\newcommand{\slopephi}{\abs{\partial\phi}}
\newcommand{\globalslopephi}{\mathfrak{L}_\lambda[\phi]}
\newcommand{\Lebae}{\mathscr{L}\text{-a.e.~}}
\newcommand{\bmu}{\boldsymbol{\mu}}
\newcommand{\bnu}{\boldsymbol{\nu}}
\let\olddownarrow\downarrow
\let\olduparrow\uparrow

\renewcommand{\downarrow}{\mathrel{\olddownarrow}}
\renewcommand{\uparrow}{\mathrel{\olduparrow}}

\usepackage{nicefrac}
 
\begin{document}

\title[On gradient descent-ascent flows in metric spaces]{On gradient descent-ascent flows in metric spaces} 

\author[N.~Isobe]{Noboru Isobe$^*$} 

\author[S.~Shimoyama]{Sho Shimoyama$^{**}$} 

\renewcommand{\thefootnote}{\fnsymbol{footnote}}
\footnote[0]{2020\textit{Mathematics Subject Classification}.
Primary 49Q22, 47H05; Secondary 49J52, 90C47, 58E30, 46T05, 60B10.}

\thanks{
$^{*}$ORCiD: \href{https://orcid.org/my-orcid?orcid=0000-0003-0512-4793}{0000-0003-0512-4793}.}
\thanks{$^{**}$Supported by {JSPS KAKENHI Grant Number JP25KJ0966}}

\address{
Addresses of N.I.
\endgraf
Center for Advanced Intelligence Project (AIP) \endgraf
RIKEN \endgraf
Meguro-ku 153-8914 \endgraf
Japan
}
\email{\href{mailto:noboru.isobe@riken.jp}{noboru.isobe@riken.jp} or \href{mailto:nobo0409@g.ecc.u-tokyo.ac.jp}{nobo0409@g.ecc.u-tokyo.ac.jp}}
\address{
Addresses of S.S.
\endgraf
Graduate School of Mathematical Sciences \endgraf
The University of Tokyo \endgraf
Meguro-ku 153-8914 \endgraf
Japan
}
\email{\href{mailto:sho-shimoyama@g.ecc.u-tokyo.ac.jp}{sho-shimoyama@g.ecc.u-tokyo.ac.jp}}

\begin{abstract}
Gradient descent-ascent (GDA) flows play a central role in finding saddle points of bivariate functionals, with applications in optimization, game theory, and robust control.
While they are well-understood in Hilbert and Banach spaces via maximal monotone operator theory, their extension to general metric spaces, particularly Wasserstein spaces, has remained largely unexplored.
In this paper, we develop a mathematical theory of GDA flows on the product of two complete metric spaces, formulating them as solutions to a system of evolution variational inequalities (EVIs) driven by a proper, closed functional $\phi$.
Under mild convex-concave and regularity assumptions on $\phi$, we prove the existence, uniqueness, and stability of the flows via a novel minimizing-maximizing movement scheme and a minimax theorem on metric spaces.
We establish a $\lambda$-contraction property, derive a quantitative error estimate for the discrete scheme, and demonstrate regularization effects analogous to classical gradient flows.
Moreover, we obtain an exponential decay bound for the Nikaidô--Isoda duality gap along the flow.
Focusing on Wasserstein spaces over Hilbert spaces, we show the global existence in time and the exponential convergence of the Wasserstein GDA flow to the unique saddle point for strongly convex-concave functionals.
Our framework unifies and extends existing analyses, offering a metric-geometric perspective on GDA dynamics in nonlinear and non-smooth settings.

\end{abstract}
\maketitle

\tableofcontents

\section{Introduction}
Gradient descent–ascent (GDA) flows are dynamical systems for finding saddle points of a bivariate functional $\phi$ on two complete metric spaces $X$ and $Y$, where $(x^\ast, y^\ast) \in \XY$ is said to be a saddle point of $\phi$ if it satisfies 
\[
    \phi(x^\ast, y) \le \phi(x^\ast, y^\ast) \le \phi(x, y^\ast)
    \text{ for any } (x, y) \in \XY.
\]
GDA flows originate from the Arrow--Hurwicz equation \cite{ArrowHurwiczUzawa58,Arrow2014} and were subsequently formalized by Rockafellar in Hilbert spaces \cite{Rockafellar71}, where maximal monotone operator theory yields well-posedness of the flow. 
More recent applications in applied mathematics---e.g.~distributionally robust optimization and reinforcement learning---have sparked interest in GDA over the space of probability measures endowed with the 2-Wasserstein metric, Wasserstein GDA.
Very recently, Wang and Chizat \cite{pmlr-v247-wang24c} posed an open problem regarding the asymptotic behavior of the Wasserstein GDA flow, specifically the convergence of the flow to a saddle point of $\phi$.

The questions raised here concern how the geometry of the ambient spaces $X$ and $Y$ influences GDA dynamics, particularly when $X$ and $Y$ are nonlinear.
The previous analysis of GDA flows lacks the geometric insight that underpins gradient flows in metric spaces. 
In the gradient-flow setting, Ambrosio--Gigli--Savaré in \cite{AGS} established the existence, uniqueness, and stability of the flows in general metric spaces, which are based on a geometric condition: the uniform convexity of the squared distance along curves.
By contrast, in the GDA flow setting, the analysis of its properties depends on the linearity and reflexivity of ambient spaces, see e.g.~\cite{Rockafellar71}, and it remains unclear what geometric properties of nonlinear spaces are necessary or sufficient to ensure analogous results.

This paper takes the first step toward filling this gap through a formulation of GDA on general metric spaces (\cref{sec:properties_of_EVIs}) and the construction of a well-posedness theory (\cref{sec:saddle_pro,sec:existence,sec:WGDA}).
In particular, the existence and stability of the GDA flow are indispensable for any convergence analysis in the Wasserstein setting.
Whereas Ambrosio--Gigli--Savaré established well-posedness of gradient flows via minimizing-movement schemes that preserve energy dissipation, in the GDA context, the key difficulty is that $\phi$ need not decrease along the flow.
Our main technical contribution is to overcome this obstacle under mild regularity on $\phi$: we develop a minimizing-maximizing movement scheme on metric spaces (\cref{sec:existence}). 
We prove that the scheme is well-defined by establishing a minimax theorem and the unique existence of saddle points in metric spaces (\cref{sec:saddle_pro}).
Finally, we apply these abstract results to Wasserstein GDA and obtain the first well-posedness theorem (\cref{sec:WGDA}).

\subsection{Problem setting}
Consider two complete metric spaces $(X,\distance_X)$ and $(Y,\distance_Y)$, and their product $Z \coloneqq \XY$ equipped with the canonical $\ell^2$ product metric $\distance_Z$.
Our focus is on a bivariate functional
\(
  \phi\colon Z\to \exR,
\)
assumed to be \emph{proper}, that is, the effective domain 
\[
    \Domain\phi\coloneqq\Set{x \in X|\phi(x, y) < +\infty \text{ for any } y \in Y} \times \Set{y \in Y | \phi(x, y) > -\infty \text{ for any } x \in X}
\]
is nonempty.
We additionally suppose that $\phi$ is \emph{closed} --\cref{def:closedness}-- and satisfies a suitable functional value condition on the complement $(\Domain \phi)^c$ --\cref{ass:function_value_on_and_out_domain}--, and the \emph{regularization $\Phi(\bullet,\bullet;x,y)$ of $\phi$ by the squared distance} is \emph{$(\lambda+{1}/{\tau})$-convex-concave} with some $\lambda\in\R$ --\cref{def:lambda_convex_concave_along_curve,assump:convex_concavity_Phi}-- for every $(x,y)\in X\times Y$ and $\tau>0$, where $\Phi$ is defined by
\begin{equation}
    \label{eq:defPhi}
    X\times Y\ni(x',y')\longmapsto\Phi_\tau(x,y;x',y')\coloneqq\phi(x',y')+\frac{1}{2\tau}\qty(\distance_X^2(x',x)-\distance_Y^2(y',y))\in\exR.
\end{equation}
Closedness and $(\lambda+{1}/{\tau})$-convex-concavity roughly require lower semicontinuity and ($\lambda$-)convexity for the first variable $x\in X$ of $\phi$, upper semicontinuity and ($\lambda$-)concavity for the second variable $y\in Y$, and the uniform or strong convexity of the squared distances $\distance^2_X$ and $\distance_Y^2$.
More precisely, it is necessary to consider cases where the effective domain is not the entire space $Z$.

For the above triplet $(X,Y,\phi)$ we introduce a formulation of {GDA flow} in metric spaces $X$ and $Y$ by the following system of variational inequalities:
\begin{restatable}[$\lambda$-evolution variational inequalities]{definition}{systemEVIs}\label{def:EVIs}
  Let $\lambda$ be a real number.
  A solution of the \emph{Evolution Variational Inequalities} $\EVSI{X}{Y}{\phi}{\lambda}$ is said to be a pair of (continuous) curves $(u, v)\colon(0,+\infty)\ni t\mapsto(u_t,v_t)\in\Domain\phi$ satisfying
  \begin{equation}
      \label{eq:EVSI}
      \left\{
      \begin{aligned}
          &\frac{1}{2} \Urd \distance_X^2(u_t,x)+\frac{\lambda}{2}\distance_X^2(u_t,x)+\phi(u_t, v_t)\leq\phi(x, v_t),\\
          &\frac12\Urd\distance_Y^2(v_t,y)+\frac{\lambda}{2}\distance_Y^2(v_t,y)+\phi(u_t,y)\leq\phi(u_t,v_t),
      \end{aligned}
      \right.
  \tag{EVIs$_\lambda$}
  \end{equation}
  for every $t\in(0,+\infty)$ and $(x,y)\in\Domain\phi$, where $\Urd$ is the right Dini derivative.
\end{restatable}
Intuitively, $u$ flows in the direction that minimizes $\phi(\bullet,v_t)$, and $v$ flows in the direction that maximizes $\phi(u_t,\bullet)$ at each time $t$.
The following example illustrates this intuition.
\begin{example}\label{ex:Hilbert}
    When $X$ and $Y$ are Hilbert, the system of inequalities is reduced to the following:
    \begin{equation}\label{eq:GDA_Hilbert}
        -\dot{u}(t) \in \partial_1 \phi(u(t),v(t)), \quad +\dot{v}(t) \in \partial_2 \phi(u(t),v(t))\quad\text{for }\Lebae t \in (0,+\infty),
    \end{equation}
    where \(\partial_1 \phi\colon X\times Y\longrightarrow2^X\) and \(\partial_2 \phi\colon X\times Y\longrightarrow 2^Y\) are the sub- and super-differentials of \(\phi\) with respect to the first and second variables, respectively.
\end{example}
\subsection{Summary of Main Results}
The main result of this paper is \cref{thm:existence}, which establishes the existence of a unique solution to $\EVSI{X}{Y}{\phi}{\lambda}$ under the following regularity assumption:
\begin{restatable}[Decomposition and regularity of \(\phi\)]{assumption}{decompAssump}\label{assump:decomp}
  There exist a functional  \(\varphi \colon X \times Y \to \R\), a lower semi-continuous functional \(\psi_X\colon X\to\overR\) and an upper semi-continuous functional \(\psi_Y\colon Y\to\underR\) such that
  \[
    \phi(x,y)=\varphi(x,y)+\psi_X(x)+\psi_Y(y),
  \]
  for any \((x, y) \in \Domain\phi\).  
  Moreover, the real-valued functional \(\varphi\) satisfies one of:
  \begin{enumerate*}[%
        label=(\roman*),               
        ref=\theassumption(\roman*)
        ]
    \item\label{assump:countinuous_coupling_mild} \(\varphi\) is \textbf{continuous}.
    \item\label{assump:countinuous_coupling} \(\varphi\) is \textbf{locally Lipschitz}.
  \end{enumerate*}
\end{restatable}
The above assumption enables the existence of solutions with initial values in the closure of the effective domain as follows:
\begin{restatable}[Well-posedness of $\EVSI{X}{Y}{\phi}{\lambda}$]{theorem}{existenceMain}\label{thm:existence}
  Suppose that Assumptions \ref{assump:countinuous_coupling}, \ref{ass:function_value_on_and_out_domain}, and \ref{assump:convex_concavity_Phi} hold for some $\lambda\in\R$.
  Then, for any $w_0=(x,y)\in\overline{\Domain\phi}$, there exists the unique solution $S[w_0] \coloneqq w=(u,v)\colon (0,+\infty) \longrightarrow X\times Y$ of $\EVSI{X}{Y}{\phi}{\lambda}$ with $\lim_{t \downarrow 0} w_t=w_0$.
  Moreover, if $ \{ w^n_0\}_{n \ge 1} \subset \overline{\Domain \phi}$ satisfies $w_0^n \to w_0 \in \overline{\Domain \phi}$ as $n \to  \infty$, then it holds that $S_t[w^n_0] \to S_t[w_0]$ as $n \to +\infty$ for any $t > 0$.
\end{restatable}
\begin{remark}
One of the main interests in this paper lies in whether the well-posedness can be established in the case of $\varphi\neq0$, which arises in Wasserstein GDA flows.
When $\varphi = 0$ i.e.~when there is no cross term between $u$ and $v$, $\EVSI{X}{Y}{\phi}{\lambda}$ reduces to two independent EVIs or gradient flows, and its well-posedness can be shown via the well-known argument using the minimizing movement as in \cite{AGS}.
However, when the cross term is nonzero, the movement in gradient flows is generally lost in the GDA flow in $\EVSI{X}{Y}{\phi}{\lambda}$.
Nevertheless, the above result asserts that a solution starting from the closure of the domain exists globally in time $t$ under the regularity condition of \cref{assump:countinuous_coupling}.
\end{remark}
\begin{remark}[Linear case]
When $\phi$ is ($0$-)convex-concave, and, $X$ or $Y$ are reflexive Banach spaces, we can employ the maximal monotone operator theory with \cite{MR285942}.
We then obtain the well-posedness result as in \cref{thm:existence} without \cref{assump:decomp}.
\end{remark}
A key technique to establish the well-posedness is constructing a discrete-time sequence $(W^n_\tau)_{n=0}^\infty\coloneqq((U^n_\tau,V^n_\tau))_{n=0}^\infty\subset X\times Y$ approximating the solution $((u_t,v_t))_{t>0}$ using the \emph{minimizing-maximizing movements scheme} formally given by as follows: for $(U^0_\tau,V^0_\tau)\coloneqq w_0$ and  $n\in\Z_{\geq0}$,
\[
(U^{n+1}_\tau,V^{n+1}_\tau)\coloneqq J_\tau(U^n_\tau,V^n_\tau)\coloneqq\argminimax_{(U',V')\in X\times Y}\Phi(U^n_\tau,V^n_\tau;U',V'),
\]
where the right-hand side is the set of all saddle points of $\Phi$.
This scheme is well-defined if $\tau\in(0,\nicefrac{1}{\lambda^-})$ with $\lambda^- \coloneqq \max\qty{-\lambda,0}$, because of the minimax theorem for saddle point variational problems on metric spaces, which is shown in \cref{sec:saddle_pro}. 
Specifically, we establish the following existence theorem for saddle points.
Notably, we do not require the compactness of $X$ and $Y$ which can be infinite-dimensional such as Wasserstein spaces.
\begin{restatable}[Unique existence of saddle point for strongly convex-concave functions]{theorem}{existenceSaddlePoint}\label{thm:existence_saddle_point}
Suppose that \cref{ass:function_value_on_and_out_domain} holds and $\phi$ is also closed and $\lambda$-convex-concave for some $\lambda > 0$.
Then, there exists the unique saddle point $(x^\ast, y^\ast) \in \Domain\phi$.
\end{restatable}

In \cref{sec:properties_of_EVIs}, we investigate structural properties of EVIs.
As a basic property, under the weak assumption of \cref{assump:countinuous_coupling_mild}, we prove the \emph{$\lambda$-contraction property} in \cref{thm:contraction:contraction}, or the stability of solutions with respect to initial values, which implies the uniqueness of the solution.
We also observe that the \emph{regularizing effect} known in gradient flows in metric spaces, as described in \cite{AGS,MURATORI2020108347}, is present in solutions to EVIs under \cref{assump:countinuous_coupling}.
Furthermore, we prove that the duality gap or the Nikaid\^o--Isoda error $\operatorname{NI}\colon\XY\to[0,+\infty]$ in \cite{NikaidoIsoda55}, which is defined by
\begin{equation}\label{eq:NikaidoIsosa}
    \operatorname{NI}(x,y)\coloneqq\sup_{(x',y')\in X\times Y}\qty(\phi(x,y')-\phi(x',y)),
\end{equation}
has an exponential bound along the EVIs-flow.
This kind of bounds are specific to the GDA setting, and  essential in algorithms for finding saddle points.

In \cref{sec:existence}, we provide a proof of the main theorem --\cref{thm:existence}--.
One of the key steps in the proof is \cref{prop:convergence}, which implies the convergence of the movements when the initial data $w_0$ belongs to a dense subset of $\overline{\Domain\phi}$.
Specifically, we derive the following error estimate of order $1/2$:
\begin{align*}
    \distance_Z^2(w_t,W_{t/n}^n)
    \le \slopephi^2(w_0)C(T, N, \lambda^-)\frac{1}{n} \text{ for any } t \in [0, T] \text{ and } n \geq N,
\end{align*}
where $T > 0 $, $N \geq 1$ such that $T/N \in \qty(0, \nicefrac{1}{\lambda^-})$, and $C(T, N, \lambda^-)$ is a constant given in \cref{prop:convergence}.
The functional $\slopephi$ will be introduced in \eqref{eq:def_slope} and is called the \emph{modified} (local) slope in this paper, which is different from the (local) slope used in e.g.~\cite{DeGiorgi1980,AGS,MURATORI2020108347}.
While the original slope represents a rate of decrease of the value of $\phi$ itself, our modified slope is distinctive in that it represents a rate of decrease of the duality gap $\operatorname{NI}$ of $\phi$.

We apply the abstract theory of \cref{thm:existence,thm:contraction} to the well-posedness and long-time behavior of the Wasserstein GDA flow.
We show that the Wasserstein GDA, defined canonically using Wasserstein (sub)differential calculus, is equivalent to EVIs in \cref{thm:equivalence}. 
Combining this equivalence with the abstract theory yields the first unique existence results for the Wasserstein GDA.
We also derive the exponential convergence of the Wasserstein GDA to saddle points for strongly convex-concave functionals; see \cref{cor:existence_of_wgda}.
This result addresses the open problem raised by \cite{pmlr-v247-wang24c}, subject to the assumptions of Hilbert structure of the base spaces of probability measures and convex-concavity of the functional $\phi$.

\subsection{Related works}\label{sec:motivation}
Our research is relevant to gradient flows on metric spaces, GDA algorithms for saddle point problems, and their applications to deep learning.
In the following, we review the previous research on each of these topics.
\paragraph{Gradient flow on metric spaces}
Clearly, when \(\phi\) in \eqref{eq:GDA_Hilbert} does not depend on the second variable i.e.~\(\phi(x,y) = \phiX(x)\) for $(x,y)\in X\times Y$, this system reduces to the well-known gradient flow given by \( -\dot{u}(t) \in \partial \phiX(u(t)) \).
One of the most characteristic properties of gradient flows is that the value of the function \(\phiX(u(t))\) decreases over time $t$ until the solution $u_t$ reaches a critical point of \(\phiX\).
Using this property, gradient flows have been employed as an algorithm for constructing minimizing sequences for minimization problems.
Notably, Ambrosio--Gigli--Savaré \cite{AGS} formulated gradient flows in metric spaces, generalizing the theory from Banach spaces by Crandall--Liggett \cite{MR287357}, and developed an abstract framework based on this formulation.
This theory has been instrumental in guaranteeing the existence of gradient flows in the Wasserstein space, as demonstrated by Jordan--Kinderlehrer--Otto~\cite{JKO98}.
It has also been applied to cases where $\phiX$ is time-dependent and to systems of two gradient flows, see \cite{Kopfer2017,Ferreira2018,MIMURA20171477,Blanchet+15}.

On the other hand, the GDA, which is our target, has a structure that is difficult to handle in these frameworks.
Specifically, the generator of GDA flow tends to increase the value of $\phi$ in the space $Y$.
As a result, the trajectory of the GDA for a convex function can be periodic: consider \eqref{eq:GDA_Hilbert} for the case $X=Y=\R$ and $\phi(x,y)=xy$.
The appearance of such non-trivial periodicity is a difficulty that does not arise in gradient flows for convex functions.

\paragraph{GDA flows for saddle point problems}\label{subsec:GDA_saddle_point}
The original form of GDA was introduced by Arrow--Hurwicz--Uzawa \cite{ArrowHurwiczUzawa58} as a sequential algorithm for approximately solving saddle point problems for bivariate functions on Euclidean spaces, and research has continued in the field of numerical analysis, optimization and evolutionary game theory; see, for example, \cite{Benzi_Golub_Liesen_2005,Arrow2014,doi:10.1080/02331934.2023.2215799,CHERUKURI201610,9483346,Lu2022osrresolutionodeframeworkunderstanding,Ryan23,Awi24}.
Here, the \emph{saddle point problem for a bivariable function $\phi$ on $X\times Y$} is given in the following form:
\[
\text{Find } (x^\ast, y^\ast) \in X \times Y \text{ such that } \phi(x^\ast, y) \le \phi(x^\ast, y^\ast) \le \phi(x, y^\ast) \text{ for all } (x, y) \in X \times Y.
\]

In general, the Convex Analysis framework developed by Rockafellar is also helpful for saddle point problems in Banach spaces, and it leads to the maximal monotonicity of the operator that generates the GDA flow in~\cite{MR285942,Rockafellar71,Rockafellar76,Asakawa86,GOSSEZ1972220}.
This shows that the theory of nonlinear semigroup on Banach spaces using maximal monotone operators (see, e.g.,~\cite{Brezis73,MR287357}) also works effectively for GDA.

Our research gives a generalization of GDA that can be considered not only in linear spaces as shown above, but also in Wasserstein spaces and Riemannian manifolds, as explained below.
\paragraph{Saddle point problems on metric spaces in Deep Learning}\label{subsec:W-GDA}
Recently, in machine learning, particularly in deep learning theory, saddle problems on metric spaces have been considered. For example, distributionally robust optimization~\cite{BAJGIRAN2022111608}, strategic distributional shift~\cite{Conger23,conger2024coupledwassersteingradientflows}, Riemannian minimax problems~\cite{Peiyuan23}, Generative Adversarial Networks (GAN)~\cite{pmlr-v97-hsieh19b} and mean-field two-player zero-sum games~\cite{NEURIPS2020_e97c864e,Conforti23} involve saddle point problems on Wasserstein spaces over Euclidean spaces or Riemannian manifolds. It is also known that the variance around the Wasserstein barycenter can be evaluated using the minimax theorem for saddle point problems~\cite{pass2023generalized}.

For mean-field two-player zero-sum games, there has been active research on GDA flows that can converge to the saddle point~\cite{NEURIPS2020_e97c864e, pmlr-v238-dvurechensky24a,kim2024symmetric,cai2024convergenceminmaxlangevindynamics,pmlr-v247-wang24c}.
The main interest here is the convergence of GDAs for the \emph{entropy-regularized functional}, which is typically given by the functional $\mathcal{L}_\beta$ on the set of absolutely continuous probability measures $\Pac(\R^{d_1})\times\Pac(\R^{d_2})$ on Euclidean spaces with the form
\begin{equation}
    \label{eq:ent_reg_objective}
    \mathcal{L}_\beta(\mu, \nu)\coloneqq\iint\limits_{\R^{d_1}\times\R^{d_2}}\ell\dd{(\mu\otimes\nu)}+\frac{1}{\beta}\qty(\relEnt\given{\mu}{\rho_{\X}}-\relEnt\given{\nu}{\rho_{\Y}}),
\end{equation}
where $\ell$ is a function on $\R^{d_1}\times\R^{d_2}$, $\beta>0$ is the so-called inverse temperature, $\relEnt$ is the relative entropy defined in \cite[Definition 9.4.1]{AGS}, and $\rho_{\X}\in\Ptwo(\R^{d_1]})$, $\rho_{\Y}\in\Ptwo(\R^{d_2})$ are some Gaussian distributions.
The standard setting in \cite{kim2024symmetric,cai2024convergenceminmaxlangevindynamics} assumes that $\ell$ is continuously differentiable, convex-concave, and that the gradient $\nabla\ell$ is Lipschitz continuous.
Note that the relative entropy \(\relEnt\) is strongly convex with respect to the first variable by the Bakry--Émery criterion~\cite{MR772092,STURM2005149}, which allows us to apply the results presented in this chapter.

Specifically, \cref{thm:existence,thm:equivalence,thm:contraction} provide the existence of Wasserstein GDA for $\mathcal{L}_\beta$ and the decreasing property of the error in \eqref{eq:NikaidoIsosa}, which has implications in~\cite{pmlr-v247-wang24c,kim2024symmetric,cai2024convergenceminmaxlangevindynamics}.
Furthermore, \cref{thm:contraction} gives the convergence of Wasserstein GDA to the (unique) saddle point of $\mathcal{L}_\beta$ because the point gives the trivial curves satisfying EVIs \eqref{eq:EVSI}.
This convergence result partially addresses the open problem in \cite{pmlr-v247-wang24c}.
As future work, we need to consider the convergence of Wasserstein GDA to a saddle point for cases where $\ell$ is $\lambda$-convex-concave with $\lambda\in\R$.

\section{Saddle point problems and minimax theorem on metric spaces}\label{sec:saddle_pro}
In this section, we establish a foundation for studying the existence of gradient descent-ascent flows.
Specifically, we extend the theory for saddle point problems on Banach spaces, which has been developed in~\cite{Rockafellar71}, to a metric space setting.
Let $(X,\distance_X)$, $(Y,\distance_Y)$ be complete metric spaces.
We equip $X \times Y$ with the $\ell^2$ product distance.
Let $\phi \colon \XY \to \exR$ be proper i.e. $\Domain\phi\coloneqq\qty(\Domain_X\phi)\times\qty(\Domain_Y\phi) \neq \emptyset$, where
\begin{equation}
    \begin{aligned}
        &\Domain_X\phi\coloneqq\Set{x\in X| \phi(x,y)<+\infty\text{ for all } y\in Y},\\
    &\Domain_Y\phi\coloneqq\Set{y \in Y| \phi(x,y)>-\infty\text{ for all } x \in X}.
    \end{aligned}
    \label{eq:effective_domain_on_X_Y}
\end{equation}
The reason for defining the effective domain $\Domain\phi$ in this way will become clear in \cref{rem:enough_work_on_A_B}.
We also follow the standard extension of functions to the entire space in convex analysis for our saddle function $\phi$.
Roughly speaking, we prescribe values outside the effective domain $\Domain\phi$ to be either $+\infty$ or $-\infty$.
Specifically, we consider the following situation:
\begin{assumption}\label{ass:function_value_on_and_out_domain}
     It holds that
     \[
     \phi(x,y) =
        \begin{cases}
             +\infty & (x,y) \in (\Domain_X\phi)^c \times \Domain_Y\phi \\
             -\infty & (x,y) \in \Domain_X\phi \times (\Domain_Y \phi)^c
        \end{cases}.         
     \]
\end{assumption}

Let us note that there are technical difficulties, as discussed in \cite{GOSSEZ1972220}, in constructing a bivariate functional $\phi$ by summing two functionals.
Let $f \colon X \to (-\infty, +\infty]$, $g \colon Y \to [-\infty, +\infty)$ be proper; that is $\Domain f \coloneqq \Set{x \in X | f(x) < +\infty} \neq \emptyset$ and $\Domain g\coloneqq\Set{y\in Y|-g(y)<+\infty}\neq\emptyset$.
We have at least two approaches to define a function on $X\times Y$ by summing $f$ and $g$.
On one hand, define $\phi_1$ by
\begin{align}\label{eq:define_phi_approach1}
\phi_1(x,y) \coloneqq \begin{cases}
    f(x) + g(y) & \text{if } (x, y) \in \Domain f \times \Domain g\\
    + \infty & \text{if } x \notin \Domain f \\
    - \infty & \text{if } x \in \Domain f, y \notin \Domain g
\end{cases}.
\end{align}
On the other hand, define $\phi_2$ by
\begin{align}\label{eq:define_phi_approach2}
\phi_2(x,y) \coloneqq \begin{cases}
    f(x) + g(y) & \text{if } (x, y) \in \Domain f \times \Domain g\\
    + \infty & \text{if } x \notin \Domain f, y \in \Domain g \\
    - \infty & \text{if } y \notin \Domain g
\end{cases}.
\end{align}
Although $\phi_1 \neq \phi_2$ on the whole space $X \times Y$, they agree on $\Domain f \times \Domain g$.
Moreover $\Domain_X\phi_1 \times \Domain_Y \phi_1 = \Domain_X\phi_2 \times \Domain_Y \phi_2 = \Domain f \times \Domain g$.
In what follows, we introduce some definitions and assumptions based on this discussion.

\subsection{Basic properties of saddle points and closure operators}
\begin{definition}[Saddle points]\label{def:abstract_minimax}
    A pair $(x^\ast, y^\ast)\in \XY$ is a \emph{saddle point} of $\phi$ if it satisfies
    \begin{align*}
        \phi(x^\ast, y) \leq \phi(x^\ast, y^\ast) \leq \phi(x, y^\ast) \text{ for all } (x,y) \in \XY.
    \end{align*}
    We denote the collections of all saddle points of $\phi$ by $\argminimax_{ (x, y) \in  X\times Y}\phi(x, y)$.
\end{definition}
It is well-known that the existence of a saddle point is equivalent to the interchangeability of inf-sup: we can find a proof in \cite[Proposition 1.2 in Chapter VI]{Ekeland99}.
\begin{proposition}\label{prop:exchangability}
    A pair $(x^\ast, y^\ast) \in \XY $ is a saddle point of $\phi$ if and only if it holds that
    \begin{align*}
        \phi(x^\ast,y^\ast)
        = \adjustlimits\inf_{x\in X} \sup_{y\in Y} \phi(x, y)
        = \adjustlimits\sup_{y\in Y} \inf_{x\in X} \phi(x, y) \in \R,
    \end{align*}
    and 
    \begin{align*}
        \sup_{y\in Y} \phi(x^\ast, y) ={}&\adjustlimits\inf_{x\in X} \sup_{y\in Y} \phi(x, y),&
        \inf_{x\in X} \phi(x, y^\ast) ={}& \adjustlimits\sup_{y\in Y} \inf_{x\in X} \phi(x, y).
    \end{align*}
\end{proposition}

As a fundamental condition for a saddle point to exist, we introduce a semi-continuity of $\phi$.
In the following definition, we fix one variable to an effective domain and require the semicontinuity with respect to the other variable.
\begin{restatable}[Closedness]{definition}{closednessDef}\label{def:closedness}
    A functional $\phi \colon X\times Y \to\exR$ is \emph{closed} if the following hold:
    \begin{itemize}
        \item For each  $x \in \Domain_X\phi$, the functional $Y \ni y \mapsto \phi(x, y)\in\underR$ is upper semicontinuous.
        \item For each  $y \in \Domain_Y\phi$, the functional $X \ni x \mapsto\phi(x, y)\in\overR$ is lower semicontinuous.
    \end{itemize}
\end{restatable}

Note that we restrict a single variable to the effective domain in \cref{def:closedness}.
Without this restriction, elementary functions such as the following example would not be semicontinuous.
\begin{example}\label{ex:not_semicontinuous}
    Let $X$, $Y$ be $[0, +\infty)$ and 
    \[
    \begin{aligned}
        f(x) \coloneqq& 
        \begin{cases}
             +\infty & \text{if }x=0\\
             1/x & \text{if }x>0
        \end{cases}
       ,&
       g(y) \coloneqq& 
       \begin{cases}
            -\infty & \text{if }y=0\\
           \log y & \text{if }y>0 
       \end{cases}
       .
    \end{aligned}
    \]
    Then $\phi_1$ (resp.~$\phi_2$) defined as in~\eqref{eq:define_phi_approach1} (resp.~\eqref{eq:define_phi_approach2}) is closed in the sense of \cref{def:closedness}.
    However, the map $X \ni x \mapsto \phi_1(x, 0)$ is not lower semicontinuous at $x=0$ and the map $Y \ni y \mapsto \phi_2(0, y)$ is not upper semicontinuous at $y=0$.
\end{example}

 To deal with semicontinuity in pathological cases such as the above example, we introduce the \emph{closure operator} $\cl_X,\cl_Y$ by
\begin{align*}
    \cl_X \phi(x,y) &\coloneqq \inf \Set{\liminf_{n \to +\infty}\phi(x_n, y) | (x_n)_n\subset X\text{ such that }x_n \to x\text{ in }X},\\
    \cl_Y \phi(x,y) &\coloneqq \sup \Set{\limsup_{n \to +\infty}\phi(x, y_n) | (y_n)_n\subset Y\text{ such that }y_n \to y\text{ in }Y},
\end{align*}
for all $(x, y) \in \XY$.
In particular, choosing $x_n=x$ or $y_n=y$ shows
\(
  \cl_X\phi\le\phi\le\cl_Y\phi.
\)
One easily checks:
    \begin{itemize}
        \item For each fixed $x \in X$, the functional $Y \ni y \mapsto (\cl_Y\phi)(x, y)$ is upper semicontinuous.
        \item For each fixed $y \in Y$, the functional $X \ni x \mapsto (\cl_X\phi)(x, y)$ is lower semicontinuous.
    \end{itemize}
It is considered difficult to satisfy both semi-continuity conditions using this closure operator.
For example, $\cl_X(\cl_Y\phi)$ and $\cl_Y(\cl_X\phi)$ are generally not upper and lower semi-continuous, respectively.
\begin{example}
    Let $X, Y, f, g$ be the same as in~\cref{ex:not_semicontinuous}.
    Then $\cl_X (\cl_Y \phi_1)$ satisfies that
    \begin{align*}
        \cl_X (\cl_Y \phi_1)(0, 0) ={}& -\infty,&\cl_X (\cl_Y \phi_1)(0, y) ={}& +\infty,& \cl_X (\cl_Y \phi_1)(x, 0) ={}& - \infty.
    \end{align*}
    Thus the map $y \mapsto (\cl_X (\cl_Y\phi_1))(0, y)$ is not upper semicontinuous at $y=0$.
    Note that $0 \notin \Domain_X(\cl_X(\cl_Y\phi_1))$.
\end{example}

Under \cref{ass:function_value_on_and_out_domain}, the following remark holds: it is sufficient for the feasible pair \((x, y)\in X\times Y\) in \cref{def:abstract_minimax} to be in \(\Domain_X\phi\times\Domain_Y\phi\), which defines $\Domain\phi$.
\begin{remark}\label{rem:enough_work_on_A_B}
    Suppose that \cref{ass:function_value_on_and_out_domain} holds.
    Then a pair $(x^\ast,y^\ast) \in \XY$ is a saddle point of $\phi$ if and only if it holds that $(x^\ast,y^\ast) \in \Domain_X\phi \times \Domain_Y\phi$ and 
    \begin{align*}
        \phi(x^\ast, y) \leq \phi(x^\ast, y^\ast) \leq \phi(x, y^\ast) \text{ for all } (x,y) \in \Domain_X \phi \times \Domain_Y \phi.
    \end{align*}
    Moreover, the following hold:
    \begin{align*}
        \adjustlimits\inf_{x\in X} \sup_{y\in Y} \phi(x,y) &{}= \adjustlimits\inf_{x\in\Domain_X\phi}\sup_{y\in Y} \phi(x,y) = \adjustlimits\inf_{x\in\Domain_X\phi}\sup_{y\in\Domain_Y\phi} \phi(x,y), \\
        \adjustlimits\sup_{y\in Y} \inf_{x\in X} \phi(x,y) &{}= \adjustlimits\sup_{y\in\Domain_Y\phi}\inf_{x\in X} \phi(x,y) = \adjustlimits\sup_{y\in\Domain_Y\phi}\inf_{x\in\Domain_X\phi} \phi(x,y).
    \end{align*}
\end{remark}
From \cref{rem:enough_work_on_A_B}, we defined $\Domain_X \phi \times \Domain_Y\phi$ as the effective domain $\Domain\phi$ of $\phi$.
Furthermore, the following result is a direct consequence of \cref{rem:enough_work_on_A_B}.
\begin{corollary}
    Let $\phi_1, \phi_2 \colon \XY \to \exR$ be functionals satisfying \cref{ass:function_value_on_and_out_domain}.
    If $\Domain \phi_1 = \Domain \phi_2$ and $\phi_1 = \phi_2$ on $\Domain \phi_1$, then the saddle points of $\phi_1$ and $\phi_2$ are coincides.
\end{corollary}

If $\phi$ is closed, then furthermore
    \begin{itemize}
        \item $\cl_Y\phi(x, y) = \phi(x,y)$ for any $(x, y) \in \Domain_X\phi \times Y$,
        \item $\cl_X \phi(x,y) = \phi(x,y)$ for $(x, y) \in X \times \Domain_Y\phi$,
        \item $\Domain_X \phi \subset \Domain_X \cl_X\phi$ and $\Domain_Y \phi = \Domain_Y \cl_X\phi$,
        \item $\Domain_X\phi = \Domain_X \cl_Y\phi$ and $\Domain_Y \phi \subset \Domain_Y \cl_Y\phi$.
    \end{itemize}
These facts imply that, under \cref{ass:function_value_on_and_out_domain}, a point \((x^*,y^*)\in X\times Y\) is a saddle point of \(\phi\) if and only if it is a saddle point of \(\cl_X\phi\) (and hence of \(\cl_Y\phi\)).

\subsection{\texorpdfstring{$\lambda$}{lambda}-convex-concavity of functions on metric spaces}\label{sec:convex-concavity}
We introduce a notion of convex-concavity of $\phi\colon \XY\to\exR$, which will be instrumental in proving the uniqueness of saddle points and hence the well-posedness of the approximation scheme for GDA, as defined in \cref{sec:existence}.
This notion is an extension of the convex-concavity in linear spaces imposed in~\cite{Rockafellar71} and $\lambda$-convexity in gradient flows on metric spaces see~\cite{AGS, MURATORI2020108347}.

Recall the definition of $\lambda$-convexity of a functional $f\colon X\to\overR$: for $x_0$, $x_1 \in \Domain f \coloneqq \Set{ x \in X | f(x) < +\infty }$, let $\gamma \colon [0, 1] \to X$ be a curve from $x_0$ to $x_1$ i.e. $\gamma$ is continuous and satisfies $\gamma_0 = x_0$ and $\gamma_1 = x_1$.
The functional $f$ is \emph{$\lambda$-convex along $\gamma$ for some $\lambda\in\R$} if
\[
    f(\gamma_t)\leq(1-t)f(\gamma_0)+tf(\gamma_1)-\frac\lambda2\ t(1-t)\distance^2_X(\gamma_0, \gamma_1) \text{ for each } t \in [0, 1].
\]
Similarly, we call $g\colon Y\to\underR$ \emph{$\lambda$-concave along a curve $\sigma\colon[0,1]\to Y$ for $\lambda\in\R$} if $-g\colon Y\to\overR$ is $\lambda$-convex along $\sigma$.

\begin{restatable}[$\lambda$-convex-concave along curves]{definition}{lambdaConvConcAlongCurves}\label{def:lambda_convex_concave_along_curve}
    Let $x_0$, $x_1\in\Domain_X\phi$, $y_0, y_1 \in \Domain_Y\phi$, and let $\gamma\colon[0,1]\to X$ and $\sigma \colon [0, 1] \to Y$ be any curve from $x_0$ to $x_1$ and $y_0$ to $y_1$, respectively.
    A function $\phi\colon X\times Y \to \exR$ is said to be \emph{$\lambda$-convex-concave} for some $\lambda\in\R$ along $\gamma$ and $\sigma$ if the following two conditions are satisfied:
    \begin{enumerate}
        \item For any $y \in Y$, the map $x \mapsto \phi(x, y)$ is $\lambda$-convex along $\gamma$.
        \item For any $x \in X$, the map $y \mapsto \phi(x, y)$ is $\lambda$-concave along $\sigma$.
    \end{enumerate}
    Moreover, $\phi$ is called \emph{$\lambda$-convex-concave for some $\lambda \in \R$} if, for any $x_0$, $x_1\in \Domain_X\phi$ and $y_0$, $y_1 \in \Domain_Y \phi$ there exist curves $\gamma$ from $x_0$ to $x_1$ and $\sigma$ from $y_0$ to $y_1$ such that $\phi$ is $\lambda$-convex-concave along $\gamma$ and $\sigma$.
    Finally, $\phi$ is \emph{strongly convex-concave} if there exists a strictly positive number $\lambda > 0$ such that $\phi$ is $\lambda$-convex-concave.
\end{restatable}
\begin{example}\label{ex:lambda_convex_concave}
    \begin{enumerate}[wide=0pt,label=(\roman*),ref=(\roman*)]
        \item\label{ex:lambda_convex_concave:f+g} Let $f \colon X \to \overR$ and $g \colon Y \to\underR$ be $\lambda_1$-convex and $\lambda_2$-concave with some $\lambda_1$, $\lambda_2\in\R$, respectively.
        We also assume that $f$ and $g$ are proper.
        Then $\phi_1$ defined by \eqref{eq:define_phi_approach1} is $\lambda_{\textup{min}}$-convex-concave with $\lambda_{\textup{min}}\coloneqq\min\qty{\lambda_1,\lambda_2}$.
        Similarly, the functional $\phi_2$ defined by~\eqref{eq:define_phi_approach2} is also $\lambda_{\textup{min}}$-convex-concave.
        \item\label{ex:lambda_convex_concave:wasserstein} Let $\X$ and $\Y$ be separable Hilbert spaces, and let $\ell \colon \X\times\Y \to \R$ be $C^{1,1}$-functional with $L \coloneqq \mathrm{Lip}(D\ell)$ where $D\ell$ is the Fr\'echet derivative of $\ell$.
        Let $X$ and $Y$ be 2-Wasserstein spaces over $\X$ and $\Y$ defined in \eqref{eq:def_W2} below, respectively.
        We define the continuous functional $\phi \colon \XY \to \R$ by
        \begin{equation}
            \phi(\mu, \nu) \coloneqq \iint_{\X\times\Y} \ell(x, y) \dd{(\mu\otimes\nu)(x,y)} \text{ for } (\mu, \nu) \in \XY.\label{eq:def_int_ell}
        \end{equation}
        Then, for some (therefore all) $\lambda \leq -2L$, $\phi$ is $\lambda$-convex-concave along \emph{generalized geodesics}: see~\cref{subsec:conv_WGDA} or \cite[Definition~9.2.2]{AGS}.
    \end{enumerate}
\end{example}

\begin{remark}
    If $\phi$ is $\lambda$-convex-concave, then for any $x_0$, $x_1 \in \Domain_X\phi$ and $y_0$, $y_1\in \Domain_Y \phi$, letting $\gamma$ and $\sigma$ be curves from $x_0$ to $x_1$ and from $y_0$ to $y_1$ along which $\phi$ is $\lambda$-convex-concave, it holds that
    \[
    \gamma_t \in \Domain_X\phi\text{ and }\sigma_t \in \Domain_Y \phi\text{ for any }t \in [0, 1].
    \]
\end{remark}

\subsection{Existence of saddle points and minimax theorem}
In this section, we consider the existence of saddle points under a strong convex-concavity assumption.
The next lemma plays an important role in the existence of saddle points.
\begin{lemma}\label{lem:infimum_is_concave}
    Suppose that \cref{ass:function_value_on_and_out_domain} holds and
    $\phi$ is also closed and $\lambda$-convex-concave for some  $\lambda\in\R$.
    Then the functional $Y \ni y \mapsto \phi^Y(y) \coloneqq \inf_{x\in X} \phi(x, y) = \inf_{x \in \Domain_X\phi} \phi(x, y)$ is $\underR$-valued, upper semicontinuous, and $\lambda$-concave.
\end{lemma}
\begin{proof} 
    Since $\phi$ is proper, $\phi^Y(y)$ belongs to $[-\infty, +\infty)$ for any $y \in Y$.
    For any $y_0, y_1 \in \Domain \phi^Y  \subset \Domain_Y\phi$, let $\sigma$ be a curve from $y_0$ to $y_1$ as in \cref{def:lambda_convex_concave_along_curve}.
    For any $x \in X$, we have
    \begin{align*}
        \phi(x, \gamma_t^Y) \geq (1-t) \phi(x, y_0) + t \phi(x, y_1) + \frac{\lambda}{2}t(1-t) \distance_Y^2(y_0, y_1),
    \end{align*}
    for all $t \in [0, 1]$.
    Taking infimum with respect to $x \in X$, we can show that $\phi^Y$ is $\lambda$-concave.
    For any $x \in \Domain_X \phi$, we have
    \begin{align*}
        \phi(x, y) \geq \limsup_{y' \to y} \phi(x, y') \geq \limsup_{y' \to y} \phi^Y(y').
    \end{align*}
    Taking infimum with respect to~$x \in \Domain_X \phi$, we can get the upper semicontinuity of $\phi^Y$.
\end{proof}

We then provide a proof of \cref{thm:existence_saddle_point}, which concerns the existence result of a saddle point as studied in \cite{Rockafellar71,Peiyuan23} without compactness assumptions.
The idea of the proof is quite similar to in~\cite[Proposition~2.1]{Ekeland99}.
The difference is that the continuity of a one-parameter family of minimizers can be shown by strong convex-concavity of the function without the compactness assumption for the underlying metric spaces.
\begin{proof}[Proof of \cref{thm:existence_saddle_point}]
    First, we show the existence of saddle points.
    Let $\phi^Y$ be the same as in \cref{lem:infimum_is_concave}.
    Note that for any $y \in \Domain_Y\phi$ the function $\phi(\bullet, y) \colon X \rightarrow (-\infty, +\infty]$ is proper, lower semicontinuous and strongly convex.
    Thus for any $y \in \Domain_Y \phi$ there exists the minimizer $x^\ast(y) \in \Domain_X\phi$ of $\phi(\bullet, y)$ such that $\phi(x^\ast(y), y) = \phi^Y(y) \in \R$.
    Similarly,  since $\phi^Y$ is strongly-concave and upper semicotinuous, there exists the maximizer $y^\ast\in\Domain\phi^Y$ of $\phi^Y$ namely
    \begin{align}
        \label{eq:equality_of_sup_inf}
        \phi^Y(y^\ast) = \phi(x^\ast(y^\ast), y^\ast) = \adjustlimits\sup_{\Domain_Y\phi} \inf_{\Domain_X\phi} \phi \in \R.
    \end{align}
    In the following we show that $(x^\ast(y^\ast), y^\ast)$ is a saddle-point of $\phi$.
    Fix any $y \in \Domain_Y\phi$ and a curve $\sigma$ from $y^\ast$ and $y$ as in \cref{def:lambda_convex_concave_along_curve}.
    It is clear that
    \begin{align*}
        \phi^Y(y^\ast) \geq {}&\phi^Y(\sigma_t)\\
        = {}&{}\phi(x^\ast(\sigma_t), \sigma_t)\\
        \geq {}&{} (1-t) \phi(x^\ast(\sigma_t), y^\ast) + t \phi(x^\ast(\sigma_t), y) \\
        \geq{}&{} (1-t)\phi^Y(y^\ast) + t \phi(x^\ast(\sigma_t), y) \text{ for each } t\in[0,1].
    \end{align*}
    Therefore we obtain
    \begin{align}
        \label{eq:intermediate_lower_of_sup_inf}
        \phi^Y(y^\ast) \geq \phi(x^\ast(\sigma_t), y) \text{ for each } t \in (0, 1].
    \end{align}
    We can also show that $\lim_{t \downarrow 0} \distance(x^\ast(\sigma_t), x^\ast(\sigma_0)) = 0$ as follows: suppose that this continuity does not hold; then we can take $\delta > 0$ and $(t_n)_{n\geq1}$ such that $t_n \downarrow 0$ and $\distance(x^\ast(\sigma_{t_n}), x^\ast(\sigma_0)) \geq \delta$ for sufficiently large $n$.
    Since, for each $n\geq1$, the function $X \ni x \mapsto \phi(x, \sigma_{t_n})$ is $\lambda$-convex and $x^\ast(\sigma_{t_n})$ belongs to $\Domain_X \phi$, we get
    \begin{align*}
        \phi(x^\ast(\sigma_{t_n}),\sigma_{t_n})\leq\frac{1}{2}\phi(x^\ast(\sigma_{t_n}),\sigma_{t_n})+\frac12\phi(x^\ast(\sigma_0),\sigma_{t_n})-\frac{\lambda}{8}\delta^2,
    \end{align*}
    namely,
    \begin{align*}
        \phi(x^\ast(\sigma_0), \sigma_{t_n}) \geq \phi(x^\ast(\sigma_{t_n}), \sigma_{t_n}) + \frac{\lambda \delta^2}{4} = \phi^Y(\sigma_{t_n}) + \frac{\lambda \delta^2}{4},
    \end{align*}
    for sufficiently large $n$.
    Since the function $[0, 1] \ni t \mapsto \phi^Y(\sigma_t)$ is continuous at $t=0$ and $\phi(x^\ast(\sigma_0),\bullet)$ is upper semicontinuous, taking $\limsup_{n\to+\infty}$ in the above inequality yields that
    \begin{align*}
        \phi^Y(\sigma_0) \geq \phi^Y(\sigma_0) + \frac{\lambda \delta^2}{4} > \phi^Y(\sigma_0).
    \end{align*}
    This is a contradiction; this shows the continuity of the map $[0, 1] \ni t \mapsto x^\ast(\sigma_t)$ at $t=0$.
    Let $x^\ast \coloneqq x^\ast(y^\ast)$.
    Taking $\liminf_{t\downarrow0}$ in~\eqref{eq:intermediate_lower_of_sup_inf} we can see that
    \begin{align} \label{eq:lower_of_sup_inf}
        \phi^Y(y^\ast) \geq \phi(x^\ast, y) \text{ for all } y\in\Domain_Y\phi.
    \end{align}
    Since $\adjustlimits\inf_{x\in\Domain_X\phi}\sup_{y\in\Domain_Y\phi}\phi(x,y) \geq \adjustlimits\sup_{y\in\Domain_Y\phi} \inf_{x\in\Domain_X\phi}\phi(x,y)$,
    combining \eqref{eq:lower_of_sup_inf} and \eqref{eq:equality_of_sup_inf}, we obtain
    \begin{align}
        \label{eq:interchangeability_of_sup_inf_and_inf_sup}
        \adjustlimits\inf_{x\in\Domain_X \phi~} \sup_{y\in\Domain_Y \phi} \phi(x,y) = \phi^Y(y^\ast) = \phi(x^\ast, y^\ast) = \adjustlimits\sup_{y\in\Domain_Y \phi~} \inf_{x\in\Domain_X \phi} \phi(x,y) \in \R.
    \end{align}
    Letting $\phi^X(x) \coloneqq \adjustlimits\sup_{y\in\Domain_Y\phi} \phi(x, y)$, we have
    \begin{align*}
        \phi^Y(y^\ast) \geq \phi^X(x^\ast) \geq \inf_{\Domain_X\phi} \phi^X = \phi^Y(y^\ast). 
    \end{align*}
    This implies that $x^\ast$ attains a minimum of $\phi^X$ i.e.
    \begin{align}
        \label{eq:x_bar_attains_inf}
        \sup_{y\in\Domain_Y\phi} \phi(x^\ast, y) = \adjustlimits\inf_{x\in\Domain_X \phi~} \sup_{y\in\Domain_Y \phi} \phi(x,y).
    \end{align}
    Similarly, we also obtain $y^\ast$ attains a maximum of $\phi^Y$ i.e.
    \begin{align*}
        \inf_{x\in\Domain_X\phi} \phi(x, y^\ast) = \adjustlimits\sup_{y\in\Domain_Y \phi~} \inf_{x\in\Domain_X \phi} \phi(x,y).
    \end{align*}
    From \cref{prop:exchangability} and \cref{rem:enough_work_on_A_B}, the pair $(x^\ast, y^\ast) \in\Domain\phi$ is a saddle point of $\phi$.
    This ends the proof of the existence of saddle points.

    Next, we show that the saddle point is unique.
    Let $(x^\ast, y^\ast)$ and $(\widetilde{x}, \widetilde{y})\in\Domain\phi$ be saddle points of $\phi$.
    Let $\gamma$ (resp.~$\sigma$) be a curve from $x^\ast$ to $\widetilde{x}$ (resp.~$y^\ast$ to $\widetilde{y}$) as in \cref{def:lambda_convex_concave_along_curve}.
    By the strong convex-concavity of $\phi$, we have
    \begin{align*}
        \frac{\lambda}{8} \distance_X^2(x^\ast, \widetilde{x}) &\leq \frac12\phi(x^\ast, y) + \frac12 \phi(\widetilde{x}, y) - \phi(\gamma_{1/2}, y), \\
        \frac{\lambda}{8} \distance_Y^2(y^\ast, \widetilde{y}) &\leq \phi(x, \sigma_{1/2}) - \frac12\phi(x, y^\ast) - \frac12 \phi(x, \widetilde{y}),
    \end{align*}
    for any $x \in X$ and $y \in Y$.
    Inserting $x = \gamma_{1/2}$ and $y = \sigma_{1/2}$ in the above inequalities and taking sum, we get
    \begin{align*}
       \distance_X^2(x^\ast, \widetilde{x}) + \distance_Y^2(y^\ast, \widetilde{y})
        \leq\frac4\lambda\qty(\qty(\phi(x^\ast,\gamma_{1/2}^Y)-\phi(\gamma_{1/2}^X,y^\ast))+\qty(\phi(\widetilde{x},\gamma_{1/2}^Y)-\phi(\gamma_{1/2}^X,\widetilde{y}))).
    \end{align*}
    Since $(x^\ast, y^\ast)$ and $(\widetilde{x}, \widetilde{y})$ are saddle points, the right-hand side of the above is less than $0$, which implies $x^\ast = \widetilde{x}$ and $y^\ast = \widetilde{y}$.
\end{proof}

\section{A metric approach to gradient descent-ascent flows: a system of evolution variational inequalities (EVIs)}\label{sec:properties_of_EVIs}
Let $X, Y$ be complete metric spaces, $X \times Y$ be equipped with $\ell^2$-distance, and $\phi \colon X \times Y \to \exR$ be a proper and closed functional.

This section aims to characterize gradient descent-ascent flows in metric spaces and present their basic properties.
In a one-variate functional setting, gradient flows in metric spaces can be characterized through a specific inequality, called \emph{evolution variational inequality} (EVI), without the linear structure of metric spaces: see~\cite{AGS, MURATORI2020108347}.
Moreover, this characterization induces some useful properties for the analysis of asymptotic behavior of flows: for more details, see \cite[Theorem~3.5]{MURATORI2020108347}.
Based on this spirit, we extend the EVI formulation for this bivariate functional setting as follows:

\systemEVIs*

The rest of this section is devoted to presenting basic properties for the solutions of $\text{EVIs}_\lambda(X,Y,\phi)$.
The existence of flows is proved in \cref{sec:existence} under a convexity-concavity assumption for functionals and distances.
In what follows, by \textit{EVIs-flow} we mean a solution of $\text{EVIs}_\lambda(X,Y,\phi)$.
Unlike the one-variable (gradient flow) case, our EVIs-flows do \emph{not}, in general, make the map $t\mapsto\phi(u_t,v_t)$ either monotone or semicontinuous.
These failures can obstruct the establishment of the existence of the flows.
Due to these difficulties, we impose the decomposition conditions on $\phi$ as in \cref{assump:decomp}.
We remark that \cref{assump:countinuous_coupling} is stronger than \cref{assump:countinuous_coupling_mild}.
\begin{example}
    The functionals $\phi_1$ and $\phi_2$ in \cref{ex:lambda_convex_concave}~\ref{ex:lambda_convex_concave:f+g} satisfy \cref{assump:countinuous_coupling}.
    Moreover, the functional $\phi$ in \cref{ex:lambda_convex_concave}~\ref{ex:lambda_convex_concave:wasserstein} satisfies \cref{assump:countinuous_coupling}, which will be proved in \cref{prop:application_functional_is_local_lip}.
    The example introduced in \eqref{eq:ent_reg_objective} also satisfies \cref{assump:countinuous_coupling};
    specifically,  $\varphi(\mu,\nu)=\int\ell\dd{(\mu\otimes\nu)}$, $\psi_X(\mu)=\beta^{-1}\relEnt\given{\mu}{\rho_{\X}}$, and $\psi_Y(\nu)=-\beta^{-1}\relEnt\given{\nu}{\rho_{\Y}}$ for $\mu\in\Ptwo(\X)$ and  $\nu\in\Ptwo(\Y)$, where $\X = \R^{d_1}$ and $\Y = \R^{d_2}$.
    For example, $\varphi(\mu,\nu)=(\int\ell\dd{(\mu\otimes\nu)})^{1/2}$, this satisfies (i) but not (ii).
\end{example}

As the following proposition shows, we can also characterize $\EVSI{X}{Y}{\phi}{\lambda}$ using the integral form.
This form is convenient when proving the existence of EVIs-flows.
\begin{proposition}[Integral characterizations of EVIs]\label{prop:integral_characterization}
    Suppose that \cref{ass:function_value_on_and_out_domain} holds.
    Let $w_\bullet=(u_\bullet,v_\bullet) \colon (0, +\infty) \to \overline{\Domain\phi}$ be a curve.
    Consider the following two conditions:
    \begin{enumerate}
        \item $w$ is a solution of $\EVSI{X}{Y}{\phi}{\lambda}$;
        \item $w$ satisfies the following two conditions:
        \begin{itemize}
            \item for every $(x,y)\in\Domain\phi$, it holds that $\phi(u_\bullet,y)$, $\phi(x,v_\bullet)$, $\phi(u_\bullet,v_\bullet)\in \Loneloc((0,+\infty))$;
            \item for all $s$, $t\in(0,+\infty)$ with $s<t$, and $(x,y) \in \Domain \phi$, it holds that
            \begin{equation}
            \label{eq:EVIs'}
            \tag{$\textup{EVIs}_\lambda'$}
                    \left\{
            \begin{aligned}
                &\frac{1}{2}\qty(\distance_X^2(u_t,x)-\distance_X^2(u_s,x))\leq\int_s^t\qty(\phi(x,v_r)-\phi(u_r,v_r)-\frac{\lambda}{2}\distance_X^2(u_r,x))\dd{r},\\
                &\frac{1}{2}\qty(\distance_Y^2(v_t,y)-\distance_Y^2(v_s,y))\leq\int_s^t\qty(\phi(u_r,v_r)-\phi(u_r,y)-\frac{\lambda}{2}\distance_Y^2(v_r,y))\dd{r}.
            \end{aligned}
            \right.
            \end{equation}
        \end{itemize}
    \end{enumerate}
    Then the following hold:
    \begin{enumerate}[label=(\alph*),font=\upshape]
        \item if \cref{assump:countinuous_coupling_mild} holds, then (i) implies (ii);
        \item if \cref{assump:countinuous_coupling} 
        holds, then (i) is equivalent to (ii).
    \end{enumerate}
    Furthermore, in the both cases, \eqref{eq:EVIs'} is equivalent to 
    \begin{equation}
        \label{eq:EVIs'exp}
        \tag{$\textup{EVIs}''_\lambda$}
                \left\{
        \begin{aligned}
            &\frac{1}{2}\qty(\e^{\lambda(t-s)}\distance_X^2(u_t,x)-\distance_X^2(u_s,x))\leq\int_s^t\e^{\lambda (r-s)}\qty(\phi(x,v_r)-\phi(u_r,v_r))\dd{r},\\
            &\frac{1}{2}\qty(\e^{\lambda(t-s)}\distance_Y^2(v_t,y)-\distance_Y^2(v_s,y))\leq\int_s^t\e^{\lambda (r-s)}\qty(\phi(u_r,v_r)-\phi(u_r,y))\dd{r}.
        \end{aligned}
        \right.
    \end{equation}
\end{proposition}
\begin{proof}
    We give the proof that condition (ii) implies (i), where we will use \cref{assump:countinuous_coupling}; otherwise, the proof is the same as that of \cite[Theorem 3.3]{MURATORI2020108347}.
    In the same way as the proof of \cite[Theorem 3.3]{MURATORI2020108347}, the function $(0,+\infty)\times(0,+\infty)\ni(t,s)\mapsto\distance_Z^2(w_t,w_s)\coloneqq\distance_X^2(u_t,u_s)+\distance_Y^2(v_t,v_s)\in\R$ is Lebesgue measurable.
    Setting $V(h)\coloneqq(\int_{t_0}^{t_1}\distance_Z^2(w_{s+h},w_{s})\dd{s})^{1/2}$ with $h>0$ and \eqref{eq:EVIs'} show that
    \begin{equation*}
        \frac{1}{2}(V(h))^2+\frac{\lambda}{2}\int_0^h(V(r))^2\dd{r}\le\int_{t_0}^{t_1}\int_0^h\qty(\phi(u_s,v_{s+r})-\phi(u_{s+r},v_s))\dd{r\dd s},
    \end{equation*}
    for $h>0$ and $t_0$, $t_1\in(0,+\infty)$ with $t_0<t_1$.
    By \cref{assump:countinuous_coupling}, the right-hand side of the above is further bounded as
    \begin{align*}
        &\int_{t_0}^{t_1}\int_0^h\qty(\phi(u_s,v_{s+r})-\phi(u_{s+r},v_s))\dd{r\dd s}\\
        ={}&\int_0^h\int_{t_0}^{t_1}\qty(\qty(\phiX(u_s)-\phiX(u_{s+r}))+\qty(\phiY(v_{s+r})-\phiY(v_s)))\dd{s\dd r}\\
        &+\int_{t_0}^{t_1}\int_0^h\qty(\contiphi(u_s,v_{s+r})-\contiphi(u_{s+r},v_s))\dd{r\dd s}\\
        \le{}&h^2W(h)+L_{t_0,t_1}\int_{t_0}^{t_1}\int_0^h\distance_Z(w_{s+r},w_s)\dd{r\dd s}\le h^2W(h)+L_{t_0,t_1}(t_1-t_0)^{1/2}\int_0^hV(r)\dd {r},
    \end{align*}
    where $W(h)\coloneqq\int_0^1\qty(\qty(\phiX\qty(x_{t_0+h\xi})-\phiX\qty(x_{t_1+h\xi}))+\qty(\phiY\qty(y_{t_1+h\xi})-\phiY\qty(y_{t_0+h\xi})))\qty(1-\xi)\dd{\xi}$ and $L_{t_0,t_1}$ is a Lipschitz constant depending on $w_t\in X\times Y$, $t\in[t_0,t_1]$.
    A version of Gronwall lemma similar to \cite[Lemma 4.1.8]{AGS} yields that there exists $C_{t_0,t_1}>0$ such that $\int_{t_0}^{t_1-h}\distance_Z^2(w_{s+h},w_s)\dd{s}\le C_{t_0,t_1}h^2$ for $h\in(0,t_1-t_0)$, then $w\in C((0,+\infty);X\times Y)$.
    By \cite[Lemma A.1]{MURATORI2020108347} and \cite[Lemma 1.3]{clement2010introduction}, $w$ also satisfies \eqref{eq:EVSI}.
\end{proof}

Similar to gradient flows, EVIs-flows exhibit some analogous results collected in~\cite[Theorem~3.5]{MURATORI2020108347,}.
In particular, the monotonicity of slope along gradient flow is generalized in EVIs-flow to the monotonicity of the following modified \emph{local slope $ \abs{\partial \phi}$ and global slope $\globalslopephi$ of $\phi$}:
\begin{align}
    \abs{\partial \phi}(x,y) & \coloneqq 
        \begin{cases}\displaystyle
            \limsup_{\substack{\Domain\phi\ni(x', y')\eqqcolon z' \to z}} \dfrac{(\phi(x,y')-\phi(x',y))^{+}}{\distance_Z(z',z)} & \text { if } z \coloneqq (x,y) \in \Domain\phi \\
            0 & \text { if } (x, y) \in \Domain\phi \text { is isolated} \\
            +\infty & \text { otherwise }
        \end{cases}, \label{eq:def_slope}\\
    \globalslopephi(x,y) & \coloneqq
        \begin{cases}\displaystyle
            \sup_{\substack{\Domain\phi\ni(x', y')\eqqcolon z' \neq z}}\qty(\dfrac{\phi(x,y')-\phi(x',y)}{\distance_Z(z',z)}+\frac\lambda2\distance_Z(z',z))^+ & \text{if } z \coloneqq (x,y) \in \Domain \phi \\
            +\infty & \text{otherwise}
        \end{cases}\nonumber,
\end{align}
where $(\bullet)^+\coloneqq\max\qty{\bullet,0}$ for $\bullet\in\exR$.

\begin{theorem}[Properties of solutions of $\EVSI{X}{Y}{\phi}{\lambda}$]
    \label{thm:contraction}
    Let $w_t \coloneqq (u_t, v_t)$, $w_t^0 \coloneqq (u^0_t,v^0_t)$ and $w_t^1 \coloneqq(u^1_t,v^1_t)$ be solutions of $\EVSI{X}{Y}{\phi}{\lambda}$. 
    Under Assumptions \ref{ass:function_value_on_and_out_domain} and \ref{assump:countinuous_coupling_mild}, the following holds:
    \begin{enumerate}[leftmargin=0mm]
        \item[] \phantomsection
      \label{thm:contraction:contraction} (\textbf{$\lambda$-contraction and uniqueness})
        \begin{align*}
            \distance_Z^2(w_t^0, w_t^1) \le \e^{-2\lambda(t-s)}\distance_Z^2(w_s^0, w_s^1) \text{ for every } 0 \le s < t < +\infty.
        \end{align*}
            In particular, for any $w_0 \in \overline{\Domain\phi}$, there exists at most one solution $w_t$ of $\EVSI{X}{Y}{\phi}{\lambda}$ with $\lim_{t\downarrow0} w_t = w_0$.
    \end{enumerate}
    Moreover, if \cref{assump:countinuous_coupling} also holds, then the following statements hold:
    \begin{enumerate}[leftmargin=0mm]
        \item[] \phantomsection
      \label{thm:contraction:right_limit} (\textbf{right limits})
            for any $t > 0$ the right derivative
            \begin{align*}
                |w_t'|_+ \coloneqq \lim_{s \downarrow t}\frac{\distance_Z(w_t, w_s)}{|t-s|}
            \end{align*}
            exists and it is finite.
            Moreover, the following identity holds:
            \begin{align}
                \label{eq:right_derivative_identity}
                |w_t'|_+ = \slopephi(u_t, v_t) = \globalslopephi(u_t, v_t) \text{ for every } t > 0.
            \end{align}

        \item[] \phantomsection
      \label{thm:contraction:regularize} (\textbf{regularizing effects})
            The map $(0, +\infty) \ni t \mapsto \phi(u_t, v_t)\in\R$ and $w$ are locally Lipschitz in $(0, +\infty)$, and the map $[0, +\infty) \ni t \mapsto \e^{\lambda t} \slopephi (u_t, v_t)\in\R$ is non-increasing and right continuous.

        \item[] \phantomsection
      \label{thm:contraction:asym} (\textbf{asymptotic behavior as $t \to +\infty$})
        if $\lambda$ is strictly positive, the duality gap $\operatorname{NI}$ defined in \eqref{eq:NikaidoIsosa} satisfies the exponential bound given by
        \begin{equation}
            \operatorname{NI}(u_t,v_t)\coloneqq\sup_{(x,y)\in X\times Y}\qty(\phi(u_t,y)-\phi(x,v_t))\le\frac{\e^{-2\lambda (t-s)}}{2\lambda}\slopephi^2(w_s) \text{ for any } 0 \le s < t .
            \label{eq:NI_decrese}
        \end{equation}
    \end{enumerate}
\end{theorem}
\begin{proof}
    We prove each statement in turn.

    (\hyperref[thm:contraction:contraction]{\textbf{$\lambda$-contraction}}) 
    The proof is the same approach as \cite[Proposition 1.1]{clement2010introduction}; however, to clarify the difference between EVI and EVIs, we present it.
    Let $J\subset(0,+\infty)$ be an open interval and set $w^0_\tau\coloneqq(u^0_t,v^0_t)$, $w^1\coloneqq(u^1_t,v_t^1)\in Z= X\times Y$ and $\varDelta_t\coloneqq\frac{\e^{2\lambda t}}{2}\distance_Z^2(w^0_t,w^1_t)$ for $t>0$.
    It is sufficient to prove that
    \[
        \limsup_{\substack{(0,\delta)\ni h\downarrow0}} \int_a^b \frac{\qty(\varDelta_{t+h}-\varDelta_t)}{h}\eta_t\dd{t}\leq0,
    \]
    for any nonnegative function $\eta\in\Ccinf(J)$, where $a$, $b$ and $\delta$ are positive numbers satisfying $\supp\eta\subset[a,b]\subset[a-\delta,b+\delta]$.
    By \cref{prop:integral_characterization}, we get that
    \begin{align}
        &
        \begin{aligned}
            \frac12\qty(\e^{\lambda (t+h)}\distance^2_Z(w^0_{t+h},z)-\e^{\lambda t}\distance^2_Z(w^0_t,z))
            \leq\int_t^{t+h}\e^{\lambda r}\qty(\phi(x,v_r^0)-\phi(u_r^0,y))\dd{r},
        \end{aligned}
        \label{eq:w^0_expEVI}\\
        &
        \begin{aligned}
            \frac12\qty(\e^{\lambda (t+h)}\distance^2_Z(w^1_{t+h},z')-\e^{\lambda t}\distance^2_Z(w^1_t,z'))
            \leq\int_t^{t+h}\e^{\lambda r}\qty(\phi(x',v_r^1)-\phi(u_r^1,y'))\dd{r},
        \end{aligned}
        \label{eq:w^1_expEVI}
    \end{align}
    for any $z=(x,y)$, $z'=(x',y')\in\Domain\phi$ and $h\in(0,\delta)$.
    From \eqref{eq:w^0_expEVI} and \eqref{eq:w^1_expEVI} for $z=w^1_{t+h}$ and $z'=w^0_t$, we obtain
    \begin{align*}
                    \int_a^b\frac{\varDelta_{t+h}-\varDelta_t}{h}\eta_t\dd{t}
            ={}&\int_a^b\frac{\e^{\lambda (t+h)}}{2h}\qty(\e^{\lambda(t+h)}\distance^2_Z(w^0_{t+h},w^1_{t+h})-\e^{\lambda t}\distance^2_Z(w^0_{t},w^1_{t+h}))\eta_t\dd{t}\nonumber\\
            &+\int_a^b\frac{\e^{\lambda t}}{2h}\qty(\e^{\lambda(t+h)}\distance^2_Z(w^0_{t},w^1_{t+h})-\e^{\lambda t}\distance^2_Z(w^0_{t},w^1_{t}))\eta_t\dd{t}\nonumber\\
            \le{}&\int_a^b\e^{\lambda (t+h)}\qty(\frac{1}{h}\int_t^{t+h}\e^{\lambda r}\qty(\phi(u^1_{t+h},v_r^0)-\phi(u_r^0,v^1_{t+h}))\dd{r})\eta_t\dd{t}\\
            &+\int_a^b\e^{\lambda t}\qty(\frac1h\int_t^{t+h}\e^{\lambda r}\qty(\phi(u^0_t,v_r^1)-\phi(u_r^1,v_t^0))\dd{r})\eta_t\dd{t}.
    \end{align*}
    Finally, the continuity of the shift operator $(0,\delta)\ni h\mapsto \phi(u^1_{\bullet+h},v_r^0)$, $\phi(u^0_{r},v_{\bullet+h}^0)\in L^1(a,b)$ (see e.g. \cite[Lemma~4.3]{Brezis2011}) and \cref{assump:countinuous_coupling_mild} result in
    \[
    \begin{aligned}
        \limsup_{h\downarrow0}\int_a^b\frac{\varDelta_{t+h}-\varDelta_t}{h}\eta_t\dd{t}
        \leq\int_a^b\e^{\lambda t}\qty(\phi(u^1_t,v^0_t)-\phi(u^0_t,v_t^1)+\phi(u_t^0,v_t^1)-\phi(u_t^1,v_t^0))\eta_t\dd{t}=0,
    \end{aligned}
    \]
    which concludes the proof of (\hyperref[thm:contraction:contraction]{\textbf{$\lambda$-contraction}}).

    (\hyperref[thm:contraction:regularize]{\textbf{regularizing effects}} and \hyperref[thm:contraction:right_limit]{\textbf{right limits}})
    The same manner can prove these claims as \cite[Theorem 3.5]{MURATORI2020108347}; however, we provide the proof because of the connection to \cref{rem:discussion_about_continuous_coupling} given in the later.
    Let $\delta_+(t) \coloneqq \limsup_{h\downarrow 0}\nicefrac{\distance_Z(w_{t+h}, w_t)}{h}$ and $\delta_-(t) \coloneqq \liminf_{h\downarrow 0}\nicefrac{\distance_Z(w_{t+h}, w_t)}{h}$.
    Noting it follows that $w_t$ is locally Lipshitz on $(0, +\infty)$ by the similar argument as in \cref{prop:integral_characterization};
    thus for $\mathscr{L}$-a.e.~$t >0$ 
    \begin{align*}
        \delta_+(t) = \delta_-(t) = |w_t'| < +\infty.
    \end{align*}
    In addition, $\lambda$-contraction yields that the map $(0, +\infty) \ni t \mapsto \e^{\lambda t}\delta_+(t)$ is non-increasing; then 
    \begin{align*}
        \delta_+(t) < +\infty \text{ for every } t > 0.
    \end{align*}
    Let $t > 0$.
    In \eqref{eq:EVIs'}, summing up two inequalities, dividing by ${t - s}$ and then taking $\limsup_{s\uparrow t}$, we have
    \begin{align*}
        \phi(u_t, y) - \phi(x, v_t) + \frac{\lambda}{2}\distance_Z^2(w_t, w) \le \delta_-(t) \distance_Z(w_t, w),
    \end{align*}
    for any $(x, y) \in \Domain(\phi)$, where we use \cref{assump:countinuous_coupling}.
    Thus, for each $t > 0$ we get
    \begin{align}
        \label{eq:regularlizing_effect_1}
        \slopephi(w_t, v_t) \le \globalslopephi(u_t, v_t) \le \delta_-(t) \le \delta_+(t).
    \end{align}
    On the other hand, in \eqref{eq:EVIs'} letting $x \coloneqq u_t$, $y\coloneqq v_t$ we have
    \begin{align}
        \label{eq:regularlizing_effect_2}
        \frac{1}{2} \frac{\distance_Z^2(w_{t+h}, w_t)}{h^2}
        \leq
        \frac{1}{h} \int_0^1 \qty(\phi(u_t, v_{t + \rho h})
        - \phi(u_{t+\rho h}, v_t)) \dd{\rho} 
        - \frac{\lambda}{2} \frac{1}{h} \int_0^1 
        \distance_Z^2(w_{t+\rho h}, w_t) \dd{\rho}.
    \end{align}
    Since the last two terms in the right-hand side vanish as $h \downarrow 0$, we only need to consider the first two terms.
    From the following simple calculations
    \begin{align*}
        \frac{\phi(u_t, v_{t+\rho h}) - \phi(u_{t + \rho h}, v_t)}{h} 
        &=
        \frac{\phi(u_t, v_{t+\rho h}) - \phi(u_{t + \rho h}, v_t)}{\distance_Z(w_{t+\rho h}, w_t)} \cdot \frac{\distance_Z(w_{t+\rho h}, w_t)}{\rho h} \rho,
    \end{align*}
    taking $\limsup_{h\downarrow0}$ in \eqref{eq:regularlizing_effect_2} yields
    \(
      \delta_+^2(t) \le \slopephi(u_t, v_t)\delta_+(t).
    \)
    Combining this with \eqref{eq:regularlizing_effect_1} we get \eqref{eq:right_derivative_identity}.

    The monotonicity of the map $t \mapsto \e^{\lambda t}\slopephi(u_t, v_t) $ follows from \eqref{eq:right_derivative_identity} and the monotonicity of the map $t \mapsto \e^{\lambda t}\delta_+(t)$.
    The local Lipschitz continuity of $t \mapsto \phi(u_t, v_t)$ is the direct consequence of the monotonicity of the map $t \mapsto \e^{\lambda t}\slopephi(u_t, v_t)$ and \eqref{eq:right_derivative_identity}.

    (\hyperref[thm:contraction:asym]{\textbf{asymptotic behavior as $t \to +\infty$}})
    By \eqref{eq:right_derivative_identity}, we see that
    \[
        \phi(u_t,y)-\phi(x,v_t)\le\distance_Z(w_t,z)\slopephi(w_t)-\frac{\lambda}{2}\distance_Z^2(w_t,z)\le\frac{1}{2\lambda}\slopephi^2(w_t),
    \]
    for all $(x,y)\in X\times Y$.
    Combining the above and the monotonicity of $t\mapsto\e^{\lambda t}\slopephi(w_t)$ gives the conclusion.
\end{proof}
We can easily check that for a saddle-point $w_0$ the map $(0, +\infty) \ni t \mapsto w_t \coloneqq w_0$ is the EVIs with $\lim_{t \downarrow 0} w_t = w_0$.
Combining this fact with $\lambda$-contraction property, we immediately get the following result:
\begin{corollary}[exponential convergence to a saddle-point]
    Suppose that Assumptions \ref{ass:function_value_on_and_out_domain} and \ref{assump:countinuous_coupling_mild} hold.
    Let $w_t$ be the solution of $\EVSI{X}{Y}{\phi}{\lambda}$ with $\lim_{t\to 0}w_t = w_0 \in \overline{\Domain \phi}$.
    If $\phi$ admits a saddle-point $w^\ast \coloneqq (x^\ast, y^\ast)$, then it holds that
    \[
        \distance_Z^2(w_t, w^\ast) \leq \e^{-2\lambda t} \distance_z^2(w_0, w^\ast) \text{ for each } t > 0.
    \]
\end{corollary}

When $\lambda \le 0$, EVIs flows does not generally converge to saddle points in general.
In particular the assumption of strict positivity in 
(\hyperref[thm:contraction:asym]{\textbf{asymptotic behavior as $t \to +\infty$}}) cannot generally be removed, as shown in the following example.
\begin{example}[Periodic orbit when $\lambda=0$]\label{ex:periodic_orbit}
Let $X=Y=\R$ and define $\phi(x,y)=xy$ for $x$, $y\in\R$.
Then the solution $(u,v)$ of $\EVSI{\R}{\R}{\phi}{0}$ is equivalent to the ODE system $\dot{u}_t=-v_t$, $\dot{v}_t=u_t$.
With the initial condition $(u_0,v_0)=(1,0)$ one obtains the explicit solution $(u_t,v_t)=(\cos t,\sin t)$ for $t\ge0$, which is periodic and does not converge to the unique saddle point $(0, 0)$ of $\phi$.
\end{example}

\section{Existence results for EVIs-flows}\label{sec:existence}
This section aims to prove the main theorem of \cref{thm:existence}.
As in \cref{sec:properties_of_EVIs}, let $X, Y$ be complete metric spaces, $Z \coloneqq X \times Y$ be equipped with $\ell^2$-distance $\distance_Z$, and $\phi \colon X \times Y \to \exR$ be proper and closed.
In addition, we define $\Domain(\slopephi)\coloneqq\Set{z\in\Domain\phi|\slopephi(z)<+\infty}$.

Before considering our bivariate cases, we recall the theory of gradient flows for a functional $f \colon X \to \overR$, which is developed in~\cite[Chapter~4]{AGS}.
One of the fundamental methods for proving the existence of flows is based on the resolvent operator $J_\tau\colon X\ni x\mapsto J_\tau[x]\in 2^X$ formally defined by
\begin{align}
    \label{eq:resolvent_one_variable}
    J_\tau[x] \coloneqq \argmin_{x' \in X} \qty{f(x') + \frac{1}{2\tau}\distance_X^2(x', x)},
\end{align} for sufficiently small $\tau > 0$.
Let $x \in X$ be a finite point of $f$ i.e. $f(x) < +\infty$.
If for some $\lambda \in \R$ the following holds:
\begin{equation}
    \begin{aligned}
        &\text{For any } x ,x_0, x_1 \in X, \text{ there exists a curve } \gamma \colon [0, 1] \to X \text{ from } x_0 \text{ to } x_1 \text{ such that } \\
        &\text{the map } X \ni x' \mapsto f(x') + \frac{1}{2\tau} \distance_X^2(x', x) \text{ is } (\tau^{-1} + \lambda)\text{-convex along $\gamma$ for any } \tau \in \left(0, \frac{1}{\lambda^-} \right),
    \end{aligned} \label{cond:grad_flow_convex_assumption}
 \end{equation}
 where $\lambda^- \coloneqq (-\lambda)^+$ and we set $\frac{1}{\lambda^-} \coloneqq +\infty$ if $\lambda^- = 0$.
Then $J_\tau[x]$ becomes single-valued for any $x \in X$ and the map
\begin{align*}
    J_{t/n}^n[x] \coloneqq \underbrace{(J_{t/n}\circ \cdots \circ J_{t/n})}_{n \text{ times}}[x]
\end{align*}
converges to the solution $u_t$ of EVI with $\lim_{t\downarrow0}u_t = x$ for each $t \in [0, +\infty)$ as $n \to+\infty$: see for more details~\cite[Chapter~4]{AGS}.

In this paper, we extend this method to our bivariate setting.
Similarly to the condition~\eqref{cond:grad_flow_convex_assumption}, we impose the following joint convexity-concavity assumption for the triplet $(\phi, X, Y)$.

\begin{restatable}[$(\tau^{-1}+\lambda)$-convex-concavity of $\Phi$]{assumption}{convexConcPhi}\label{assump:convex_concavity_Phi}
Suppose that for any $z = (x, y) \in \XY$, $z_0=(x_0, y_0)$, and $z_1 = (x_1, y_1)\in \Domain \phi$ there exist curves $\gamma \colon [0, 1] \to X$ from $x_0$ to $x_1$ and $\sigma \colon [0, 1] \to Y$ from $y_0$ to $y_1$ such that
\begin{align*}
    (x', y') \mapsto \Phi_\tau(x,y; x', y')  \text{ is } (\tau^{-1}+\lambda)\text{-convex-concave along } \gamma \text{ and } \sigma \text{ for any } 0 < \tau < \frac{1}{\lambda^-},
\end{align*}
where $\Phi_\tau$ is defined by \eqref{eq:defPhi}.
Note that $\Domain \Phi_\tau(x,y;\bullet,\bullet) = \Domain \phi$ for $(x,y)\in\XY$.
\end{restatable}

\begin{example}
    Let us reconsider the cases in \cref{ex:lambda_convex_concave}. 
    As the following example shows, \cref{assump:convex_concavity_Phi} imposes strong convexity on the squared-distance functions $\distance_X^2$, $\distance_Y^2$ rather than on $\phi$ itself.
    \begin{enumerate}[wide=0pt]
        \item Consider \cref{ex:lambda_convex_concave}~\ref{ex:lambda_convex_concave:f+g} when $X$, $Y$ be CAT$(0)$: see~\cite{Bacak+2014}.
        In this setting, it is known that for any point $x\in X$ the map $X\ni x'\mapsto\distance_X^2(x',x)/2$ is $1$-convex along every minimal geodesic, and likewise on $Y$.
        Hence if $f$ and $g$ are themselves $\lambda$-convex (resp.\ $\lambda$-concave) along minimal geodesics, then the functional $\phi_1$ and $\phi_2$ defined by \eqref{eq:define_phi_approach1} and \eqref{eq:define_phi_approach2} satisfies the assumption with the same $\lambda$.
        \item Consider the Wasserstein case in \cref{ex:lambda_convex_concave}~\ref{ex:lambda_convex_concave:wasserstein}.
        Note that, as shown in \cite[Subsection 9.1]{AGS}, the squared $2$-Wasserstein distance is $1$-convex along generalized geodesics employed later in \cref{subsec:conv_WGDA}, but never along minimal geodesics.
        By virtue of the generalized geodesic, $\phi$ in \eqref{eq:def_int_ell} and $(X \times Y, \distance_Z)$ satisfies \cref{assump:convex_concavity_Phi} with any $\lambda \le -2L$.
    \end{enumerate}
\end{example}

\subsection{A resolvent operator via saddle point problem}\label{subsec:resolvent_via_saddle_point_prob}
In this section, for sufficiently small $\tau > 0$, we define the resolvent operator $J_\tau \colon X \times Y \to X\times Y$  corresponding to \eqref{eq:resolvent_one_variable} in our bivariate setting.
Then, we show that the map $\overline{\Domain \phi} \ni z \mapsto J_\tau z$ is continuous and satisfies a discrete version of EVIs in \cref{lem:discrete_VI}.
These results are used to show that the map
\begin{align}
    \label{eq:n_times_composition_Jtau}
    \overline{\Domain \phi} \ni z \mapsto J_{\tau}^nz\coloneqq(\underbrace{J_{\tau}\circ\cdots\circ J_{\tau}}_{n\text{ times}})z
\end{align}
converges as $n \to +\infty$ in \cref{prop:convergence}.

\begin{definition}[Resolvent operator and Moreau--Yosida regularization]\label{def:resolvent}
    Let $\Phi_\tau$ be the function introduced in \eqref{eq:defPhi}.
    The map $J_\tau \colon X\times Y \to \Domain \phi$ is defined by, for any $z \in X\times Y$, 
    \begin{align*}
        J_\tau z\coloneqq\argminimax_{ z' \in X\times Y}\Phi_\tau(z; z').
    \end{align*}
    In addition we define the map $\phi_\tau \colon X\times Y \to \R$ by
    \(
        \phi_\tau(z) \coloneqq \Phi_\tau(z; J_\tau z)
    \)
    for any $z \in X\times Y$.
     Denoted by $J_\tau^X z$ and $J_\tau^Y z$ the projection of $J_\tau z$ onto $X$ and $Y$, respectively.
\end{definition}

\begin{remark}\label{rem:resolvent_single_value}
    By virtue of \cref{thm:existence_saddle_point}, if \cref{ass:function_value_on_and_out_domain,assump:convex_concavity_Phi} hold for some $\lambda \in \R$, then $J_\tau z$ is single-valued for every $z \in X\times Y$ and $\tau \in \left(0, \frac{1}{\lambda^-}\right)$.
\end{remark}

The following lemma is prepared to establish the continuity of $J_\tau$ in \cref{lem:conti_J}.
\begin{lemma}\label{lem:boundedness_x_tau}
    Suppose that \cref{assump:convex_concavity_Phi} holds for some $\lambda \in \R$.
    Let  $z_\ast=(x_\ast,y_\ast)\in\Domain\phi$, $\tau>0$ and $\tau_\ast>\tau$ satisfy $1+\lambda\tau>0$.
    Then, for all $(x,y)$, $(x',y')\in X\times Y$, it holds that
    \begin{align*}
        &\distance_X^2(x',x)\le\frac{4\tau_\ast\tau}{\tau_\ast-\tau}\qty(\Phi_\tau^X(x,y;x')-\phi_{\tau_{\ast}}^Y(x_\ast;y_\ast)+\frac{1}{2\tau}\distance_Y^2(y,y_\ast)+\frac{1}{\tau_\ast-\tau}\distance_X^2(x,x_\ast)),\\
        &\distance_Y^2(y',y)\le\frac{4\tau_\ast\tau}{\tau_\ast-\tau}\qty(-\Phi_\tau^Y(x,y;y')+\phi_{\tau_{\ast}}^X(y_\ast;x_\ast)+\frac{1}{2\tau}\distance_X^2(x,x_\ast)+\frac{1}{\tau_\ast-\tau}\distance_Y^2(y,y_\ast)),
    \end{align*}
    where $\phi_\tau^X(y;x')$, $\Phi_\tau^X(x,y;x')$, $\phi_\tau^Y(x;y')$, and $\Phi_\tau^Y(x,y;y')$ are defined by
    \begin{equation}
        \label{eq:partial_Moreau-Yosida}
        \begin{aligned}
            &\phi_\tau^X(y;x')\coloneqq\sup_{y'\in Y}\qty{\phi(x',y')-\frac{1}{2\tau}\distance_Y^2(y',y)},&&\Phi_\tau^X(x,y;x')\coloneqq\phi_\tau^X(y;x')+\frac{1}{2\tau}\distance_X^2(x',x),\\
            &\phi_\tau^Y(x;y')\coloneqq\inf_{x'\in X}\qty{\phi(x',y')+\frac{1}{2\tau}\distance_X^2(x',x)},&&\Phi_\tau^Y(x,y;y')\coloneqq\phi_\tau^Y(x;y')-\frac{1}{2\tau}\distance_Y^2(y',y).
        \end{aligned}
    \end{equation}
\end{lemma}
\begin{proof}
    We only prove the first inequality; the second one follows from a similar argument.
    By the elemental inequality
    \begin{align*}
        (a + b)^2 \le (1+\varepsilon)a^2 + (1+\varepsilon^{-1})b^2 \text{ for any } a, b \ge0 \text{ and } \varepsilon > 0,
    \end{align*}
    we have
    \begin{align*}
        \Phi_\tau^X(x,y;x')={}&\phi_\tau^X(y;x')+\frac{\tau_\ast-\tau}{2\tau(\tau+\tau_\ast)}\distance_X^2(x',x)+\frac{1}{\tau+\tau_\ast}\distance_X^2(x',x)\\
        \ge{}&\phi_\tau^X(y;x')+\frac{\tau_\ast-\tau}{4\tau\tau_\ast}\distance_X^2(x',x)+\frac{1}{2\tau_\ast}\distance^2_X(x',x_\ast)-\frac{1}{\tau_\ast-\tau}\distance^2_X(x,x_\ast)\\
        \ge{}&\phi(x',y_\ast)-\frac{1}{2\tau}\distance_Y^2(y_\ast,y)+\frac{\tau_\ast-\tau}{4\tau\tau_\ast}\distance_X^2(x',x)+\frac{1}{2\tau_\ast}\distance^2_X(x',x_\ast)-\frac{1}{\tau_\ast-\tau}\distance^2_X(x,x_\ast)\\
        \ge{}&\phi^Y_{\tau_\ast}(x_\ast;y_\ast)-\frac{1}{2\tau}\distance_Y^2(y_\ast,y)+\frac{\tau_\ast-\tau}{4\tau\tau_\ast}\distance_X^2(x',x)-\frac{1}{\tau_\ast-\tau}\distance^2_X(x,x_\ast),
    \end{align*}
    where we use $\phi_\tau^X(y;x')\ge\phi(x',y_\ast)-\frac{1}{2\tau}\distance_Y^2(y_\ast,y)$ in the second inequality, and $\phi(x';y_\ast)+\frac{1}{2\tau_\ast}\distance_X^2(x',x_\ast)\ge\phi_{\tau_\ast}^Y(x_\ast;y_\ast) $ in the last inequality. 
\end{proof}

\begin{lemma}
    Suppose that \cref{ass:function_value_on_and_out_domain,assump:convex_concavity_Phi} hold for some $\lambda \in \R$.
    Then for any $(x', y') \in \Domain \phi$ the following hold:
    \begin{enumerate}
        \item $\phi_\tau^Y(x; y') > -\infty$ for any $x \in X$;
        \item $\phi_\tau^X(y; x') < +\infty$ for any $y \in Y$.
    \end{enumerate}
\end{lemma}
\begin{proof}
    We only prove (i); to show (ii), it is enough to apply the same argument for $-\phi$.
    Let $g \colon Y \to (-\infty, +\infty]$ to be that $g(y) \coloneqq \phi(x', y)$; note that $g$ is lower semicontinuous, proper and $D(g) = \Domain_X\phi$ because of~\cref{ass:function_value_on_and_out_domain}.
    Since \cref{assump:convex_concavity_Phi} holds for some $\lambda \in \R$, ones can check that \cite[Assumption~2.4.5]{AGS} holds for $g$ and same $\lambda$.
    We can get the desired result by \cite[Lemma~2.4.8]{AGS}.
\end{proof}

\begin{lemma}\label{lem:conti_J}
    Suppose that \cref{ass:function_value_on_and_out_domain,assump:convex_concavity_Phi} hold for some $\lambda \in \R$.
    The resolvent operator $J_\tau\colon\overline{\Domain\phi}\longrightarrow\Domain\phi$ is continuous for each $\tau\in(0,\lambda^-)$.
\end{lemma}
\begin{proof}
    We show that $J_\tau^X\colon\overline{\Domain\phi}\longrightarrow\Domain_X\phi$ is continuous; the continuity of $J^Y_\tau\colon\overline{\Domain\phi}\longrightarrow\Domain_Y\phi$ follows in the same way.
    Let $\overline{\Domain\phi} \ni z_n \coloneqq (x_n, y_n) \to z \coloneqq (x,y) \in \overline{\Domain\phi}$ as $n \to+\infty$.
    Additionally, we set $x^n_\tau\coloneqq J_\tau ^Xz^n$, $y^n_\tau\coloneqq J_\tau^Y z^n$ and $\phi_\tau(x^n,y^n) \coloneqq \Phi_\tau(x^n, y^n; x^n_\tau, y^n_\tau)$.
    Note that $((x^n_\tau,y^n_\tau))_n$ is bounded by \cref{lem:boundedness_x_tau},  $x^n_\tau \in \argmin_{x'\in X}\Phi_\tau^X(x^n,y^n;x')$ and $y^n_\tau \in \argmax_{y'\in Y}\Phi^Y_\tau(x^n,y^n;y')$ because of \cref{prop:exchangability}.
    We first prove that
    \begin{equation} \label{eq:first_for_continuity}
        \limsup_{n\to\infty}\Phi^X_\tau(x,y;x^n_\tau)
        \le \limsup_{n\to\infty}\Phi_\tau^X(x^n,y^n;x^n_\tau)
        \qty(=\limsup_{n\to\infty} \phi_\tau(x^n, y^n)).
    \end{equation}
    By the definition of $\Phi^X_\tau$, we see that
    \begin{equation}
        \label{eq:first_1}
        \begin{aligned}
            \limsup_{n\to\infty}\qty(\Phi_\tau^X(x,y;x^n_\tau)-\Phi_\tau^X(x^n,y^n;x^n_\tau))\le&\limsup_{n\to\infty}\qty(\phi_\tau^X(y;x^n_\tau)-\phi^X_\tau(y^n;x^n_\tau))\\
            &+\frac{1}{2\tau}\limsup_{n\to\infty}\qty(\distance^2_X(x^n_\tau,x)-\distance_X^2(x^n_\tau,x^n)).
        \end{aligned}
    \end{equation}
    Because of the boundedness of $(x_\tau^n)_n$, the second term on the right-hand side of \eqref{eq:first_1} is equal to $0$.
    By the definition of $y^n_\tau$, it holds that
    \begin{align*}
        \phi(x^n_\tau,y^n)-\frac{1}{2\tau}\distance_Y^2(y^n_\tau,y^n)
        = \phi_\tau^X(y^n;x^n_\tau)
        \ge \phi(x^n_\tau,\widetilde{y})-\frac{1}{2\tau}\distance_Y^2(\widetilde{y},y^n),
    \end{align*}
    for all $\widetilde{y}\in Y$.
    Note that we can take a bounded sequence $({\widetilde{y}}^n)_n\subset Y$ such that 
    \[
    \phi(x^n_\tau,{\widetilde{y}}^n)-\frac{1}{2\tau}\distance_Y^2({\widetilde{y}}^n,y)\ge\phi^X_\tau(y;x^n_\tau)-\frac{1}{n}. 
    \]
    Then we get
    \[
        \phi_\tau^X(y^n;x^n_\tau)\ge\phi_\tau^X(y;x^n_\tau)\underbrace{-\frac{1}{n}+\frac{1}{2\tau}\qty(\distance_Y^2(\widetilde{y}^n,y)-\distance_Y^2(\widetilde{y}^n,y^n))}_{{\to0\text{ as }n\to\infty}},
    \]
    which implies the first term of \eqref{eq:first_1} is less than $0$.
    We next prove that 
    \begin{equation}
        \lim_{n\to\infty}\phi_\tau(x^n,y^n)=\phi_\tau(x,y).
        \label{eq:continuity_phi_tau}
    \end{equation}
    For every $\widetilde{x}\in \Domain_X\phi$, we have
    \begin{align*}
        \limsup_{n\to\infty}\phi_\tau(x^n,y^n)&=\limsup_{n\to\infty}\Phi_\tau^X(x^n,y^n;x^n_\tau)
        \le \limsup_{n\to\infty}\Phi_\tau^X(x^n,y^n;\widetilde{x})
        =\Phi_\tau^X(x,y;\widetilde{x}),
    \end{align*}
    where we use the continuity of $y\mapsto\phi_\tau^X(y;\widetilde{x})$ for each $\widetilde{x}\in \Domain_X\phi$; see \cite[Lemma~3.1.2]{AGS}.
    Taking the infimum of the right-hand side above with respect to $\widetilde{x}\in \Domain_X\phi$ gives us the upper semicontinuity of $\phi_\tau$.
    Applying the same argument as above to $\Phi^Y_\tau$ also yields the lower semicontinuity of $\phi_\tau$.

    Combining \eqref{eq:first_for_continuity} and \eqref{eq:continuity_phi_tau} results in $\limsup_{n\to\infty}\Phi_\tau^X(x,y;x^n_\tau)\le\phi_\tau(x,y)$.
    Therefore, $(x^n_\tau)_n$ is a minimizing sequence of $\Phi_\tau^X(x,y;\bullet)$ and converges to the unique minimizer $x_\tau\coloneqq J_\tau^X z$.
\end{proof}


\begin{lemma}[Discretized EVIs]
    \label{lem:discrete_VI}
    Suppose that \cref{ass:function_value_on_and_out_domain,assump:convex_concavity_Phi} hold for some $\lambda \in \R$.
    Then, it holds that
    \[
        \begin{aligned}
            &\frac{\distance^2_X(J_\tau^Xz,x')-\distance_X^2(x,x')}{2\tau}+\frac{\lambda}{2}\distance^2_X(J_\tau^Xz,x')+\phi(J_\tau^Xz,J_\tau^Yz)\leq\phi(x',J_\tau^Yz)-\frac{1}{2\tau}\distance_X^2(J_\tau^X z,x),\\
            &\frac{\distance^2_Y(J_\tau^Yz,y')-\distance_Y^2(y,y')}{2\tau}+\frac{\lambda}{2}\distance^2_Y(J_\tau^Yz,y')+\phi(J_\tau^Xz,y')\le\phi(J_\tau^Xz,J_\tau^Yz)-\frac{1}{2\tau}\distance_Y^2(J_\tau^Yz,y),
        \end{aligned}
    \]
    for any $z=(x,y) \in \XY$ and $z'=(x',y')\in\Domain\phi$.
\end{lemma}

\begin{proof}
    For any $z \coloneqq (x,y) \in \XY$, since $J_\tau z$ is the saddle point of $\Phi_\tau(z; \cdot)$ and \cref{assump:convex_concavity_Phi} holds, we have
    \[
        \begin{aligned}
            &0\leq-\Phi_\tau(z;J_\tau ^Xz,J_\tau^Y z)+\Phi_\tau(z;x',J_\tau^Yz)-\frac{1+\lambda\tau}{2\tau}\distance^2_X(J_\tau^Xz,x'),\\
            &0\geq-\Phi_\tau(z;J_\tau^Xz,J_\tau^Yz)+\Phi_\tau(z;J_\tau^X,y')+\frac{1+\lambda\tau}{2\tau}\distance^2_Y(J_\tau^Yz,y'),
        \end{aligned}
    \]
    for all $(x',y')\in \Domain\phi$.
    This completes the proof.
\end{proof}

\subsection{Estimate of the slope \texorpdfstring{$\abs{\partial \phi}$}{slope} for the resolvent operator \texorpdfstring{$J_\tau$}{Jtau}}\label{subsec:slope_estimate_for_resolvent}

This section is devoted to derive an estimate of the slope $\abs{\partial \phi}$ for $J_\tau$, which is used to show the convergence of \eqref{eq:n_times_composition_Jtau} as $n\to +\infty$.
To this end, we start by proving the following lemma.
The proof is similar to \cite[Theorem~2.4.9]{AGS}.

\begin{lemma}[Local slope is global slope]\label{lem:global_slope}    
    Suppose that \cref{assump:convex_concavity_Phi} holds for some $\lambda \in \R$.
    Then it holds that
    $
        \abs{\partial\phi}(z)=\mathfrak{L}_\lambda[\phi](z)
    $
    for every $z\in\Domain\phi$.
    In particular, $\abs{\partial\phi}\colon X\times Y\longrightarrow[0,+\infty]$ is lower semicontinuous.
\end{lemma}
\begin{proof}
    If $z=(x,y)\in\Domain\phi$, we see that $\slopephi(z)\leq\mathfrak{L}_\lambda[\phi](z)$.
    Thus, it remains to prove the opposite direction.
    Fix any $\Domain\phi \ni z' \coloneqq (x', y') \neq z$.
    By \cref{assump:convex_concavity_Phi}, there exists a curve $\gamma=(\gamma^X,\gamma^Y)\colon[0,1]\longrightarrow X\times Y$ with $\gamma_0=z$ and  $\gamma_1=z'$ satisfying
    \begin{subequations}\label{eq:tangent_ineq}
    \begin{align}
        &\phi(\gamma_t^X,y)+\frac{1}{2\tau}\distance^2_X(\gamma_t^X,x)\le(1-t)\phi(x,y)+t\phi(x',y)+\frac{t}{2\tau}\qty(t-\lambda\tau(1-t))\distance_X^2(x',x),\label{eq:Xdirection}\\
        &\phi(x,\gamma_t^Y)-\frac{1}{2\tau}\distance^2_Y(\gamma_t^Y,y)\ge(1-t)\phi(x,y)+t\phi(x,y')-\frac{t}{2\tau}\qty(t-\lambda\tau(1-t))\distance_Y^2(y',y),\label{eq:Ydirection}
    \end{align}
    \end{subequations}
    for all $t\in[0,1]$ and all $\tau\in(0, \frac{1}{\lambda^-})$.
    Multiplying each of \eqref{eq:tangent_ineq} by $2\tau$ and taking the limit as $\tau\downarrow0$ yields
    \begin{equation}
        \distance_Z(\gamma_t,z)\leq t~\distance_Z(z',z).
        \label{eq:relaxed_geodesic}
    \end{equation}
    In addition, subtracting \eqref{eq:Xdirection} from \eqref{eq:Ydirection} gives
    \begin{equation}
        \phi(x,\gamma_t^Y)-\phi(\gamma_t^X,y)\geq t\qty(\phi(x,y')-\phi(x',y)+\frac{1}{2\tau}\qty(\lambda\tau(1-t)-t)\distance_Z^2(z',z)),
        \label{eq:slope_estimate}
    \end{equation}
    Thus, dividing both sides of \eqref{eq:slope_estimate} by \eqref{eq:relaxed_geodesic} results in
     \[
        \frac{\phi(x,\gamma_t^Y)-\phi(\gamma_t^X,y)}{\distance_Z(\gamma_t,z)}\geq \frac{\phi(x,y')-\phi(x',y)}{\distance_Z(z',z)}+\frac{1}{2\tau}\qty(\lambda\tau(1-t)-t)\distance_Z(z',z).
     \]
     Therefore, taking $\limsup_{t\downarrow 0}$ and then taking supremum with respect to $z'$, we obtain
     \[
     \slopephi(z)\geq\limsup_{t\downarrow0}\frac{\phi(x,\gamma_t^Y)-\phi(\gamma_t^X,y)}{\distance_Z(\gamma_t,z)}\ge\mathfrak{L}_\lambda[\phi](z),
     \]
     which concludes the proof.
\end{proof}

We are now in a position to present an estimate of the slope $\abs{\partial \phi}$.
\begin{lemma}[Slope estimate for convex-concave-functionals]\label{lem:sloepe_estimate}
Suppose that \cref{ass:function_value_on_and_out_domain,assump:convex_concavity_Phi} hold for some $\lambda \in \R$.
For any $z\in \XY$ and $\tau \in \left(0, \frac{1}{\lambda^-} \right)$ it holds that
\begin{align*}
    \abs{\partial\phi}(J_\tau z)
    \le \frac{\distance_Z(J_\tau z, z)}{\tau}
    \le \frac{1}{1+\lambda\tau}\abs{\partial\phi}(z).
\end{align*}
In particular, $J_\tau z \in D(\abs{\partial\phi})$.
\end{lemma}
\begin{proof}
    Let $z\in\Domain(\slopephi)$.
    \cref{def:resolvent} yields that for any $z'\in\Domain\phi$
    \begin{align*}
        \phi(J_\tau^X z,y')-\phi(x',J_\tau^Yz)\le\frac{1}{2\tau}\qty(\distance^2_Z(z,z')-\distance_Z^2(J_\tau z,z))
        \le\frac{\distance_Z(J_\tau z,z')}{2\tau}\qty(\distance_Z(z,z')+\distance_Z(J_\tau z,z)).
    \end{align*}
    Thus, we obtain 
    \begin{equation}
        \label{eq:wk_slope_estimate}
        \slopephi(J_\tau z) = \limsup_{\substack{z'\to J_\tau z,\\z'\in\Domain\phi}}\frac{( \phi(J_\tau^X z,y')-\phi(x',J_\tau^Yz))^+}{\distance_Z(J_\tau z,z')}\le\frac1\tau\distance_Z(J_\tau z,z).
    \end{equation}
    Furthermore, combining \cref{lem:global_slope,lem:discrete_VI}, we have
    \begin{equation}
    \begin{aligned}
        \frac{1+\lambda\tau}{2\tau}\distance_Z^2(J_\tau z,z)\le{}&\phi(x,J_\tau^Y z)-\phi(J_\tau^X z,y)-\frac{1}{2\tau}\distance_Z^2(J_\tau z,z)\\
        \le{}&\distance_Z(J_\tau z,z)\slopephi(z)-\frac{1+\lambda\tau}{2\tau}\distance_Z^2(J_\tau z,z)\\
        \le{}&\frac{\tau}{2(1+\lambda\tau)}\slopephi^2(z).
    \end{aligned}
    \label{eq:bound_distance_by_slope}
    \end{equation}
    Combining \eqref{eq:wk_slope_estimate} and \eqref{eq:bound_distance_by_slope} completes the proof.
\end{proof}

The following auxiliary lemma is used to extend the starting point of the solution of EVIs from $\Domain(\slopephi)$ to $\overline{\Domain \phi}$ in \cref{thm:existence}.
\begin{lemma}\label{cor:J_tau_z_to_z}
    Suppose that the same assumptions as in \cref{lem:discrete_VI} hold.
    For any $z\in\overline{\Domain\phi}$, it holds that $\lim_{\tau\downarrow0}J_\tau z=z$.
    In particular, $\overline{\Domain\phi}=\overline{\Domain(\slopephi)}$.
\end{lemma}
\begin{proof}
    Let $\tau$, $\tau_\ast \in \left(0, \nicefrac{1}{\lambda^-}\right)$ be satisfying $\tau_\ast  > \tau$.
    By \cref{lem:boundedness_x_tau}, we obtain
    \begin{align*}
        &\distance_X^2(J_\tau^Xz,x)\le\frac{4\tau_\ast\tau}{\tau_\ast-\tau}\qty(\Phi^X_\tau(x,y;J_\tau^Xz)-\phi_{\tau_\ast}^Y(x';y')+\frac{1}{2\tau}\distance_Y^2(y,y')+\frac{1}{\tau_\ast-\tau}\distance_X^2(x,x')),\\
        &\distance_Y^2(J_\tau^Yz,y)\le\frac{4\tau_\ast\tau}{\tau_\ast-\tau}\qty(-\Phi_\tau^Y(x,y;J_\tau^Yz)+\phi_{\tau_{\ast}}^X(y';x')+\frac{1}{2\tau}\distance_X^2(x,x')+\frac{1}{\tau_\ast-\tau}\distance_Y^2(y,y')).
    \end{align*}
    for all $z'\in\Domain\phi$.
    Note that \cref{prop:exchangability} yields that $\Phi^X_\tau(x,y;J_\tau^Xz)=\Phi_\tau^Y(x,y;J_\tau^Yz)$, then
    \begin{align*}
        \limsup_{\tau\downarrow0}\distance_Z^2(J_\tau z,z)\le&\limsup_{\tau\downarrow0}\frac{4\tau_\ast\tau}{\tau_\ast-\tau}\qty(\phi_{\tau_\ast}^X(y';x')-\phi_{\tau_\ast}^Y(x';y')+\frac{\tau_\ast+\tau}{2\tau(\tau_\ast-\tau)}\distance_Z^2(z,z'))
        =2\distance_Z^2(z,z'),
    \end{align*}
    for all $z'\in\Domain\phi$.
    Because $z\in\overline{\Domain\phi}$, we obtain $\lim_{\tau\downarrow0}\distance_Z(J_\tau z,z)=0$.
\end{proof}

\subsection{Existence of the solution of EVIs}
In this section, we give the proof of the main theorem, \cref{thm:existence}.
For the convenience of readers, we restate this:
\existenceMain*
Note that the continuity of solutions with respect to initial points directly follows from $\lambda$-contraction in  \cref{thm:contraction:contraction}.
Thus we need only to show the existence.
We first establish the existence of solutions for $w_0 \in \Domain(\slopephi)$, and then extend this result for the case $w_0 \in \overline{\Domain \phi}$.
We note that \cref{thm:existence} is the first point in \cref{sec:existence} where the technical assumption \cref{assump:countinuous_coupling} is imposed; in other words, the results proved in \cref{subsec:resolvent_via_saddle_point_prob,subsec:slope_estimate_for_resolvent} hold without this assumption.
Moreover, if the start point $w_0$ belongs to $\Domain(\slopephi)$, we remark that the existence of the solution holds under the mild assumption: \cref{assump:countinuous_coupling_mild}.

To prove the existence of solutions, for any $n \in \N$ and $\tau \in \qty(0, \frac{1}{\lambda^-})$ we define the map $J_\tau^n \colon \overline{\Domain \phi} \to \Domain\phi$ to be that
\begin{align*}
    J_{\tau}^nz\coloneqq(\underbrace{J_{\tau}\circ\cdots\circ J_{\tau}}_{n\text{ times}})z \text{ for each } z \in \overline{\Domain \phi}.
\end{align*}
\cref{rem:resolvent_single_value} implies that this map is well-defined under \cref{ass:function_value_on_and_out_domain,assump:convex_concavity_Phi}.

\begin{lemma}\label{lem:distance_between_resolvents}
Suppose that \cref{ass:function_value_on_and_out_domain,assump:convex_concavity_Phi} hold for some $\lambda \in \R$.
 Then for any $z \in \Domain(\slopephi)$, $\tau_0, \tau_1 \in \left(0, \frac{1}{\lambda^-} \right)$ and $n, m \in \N$ it holds that
 \begin{align*}
     \distance_Z^2(J_{\tau_0}^n z, J_{\tau_1}^m z) 
     \leq{}& \slopephi^2(z) C(\tau_0, \tau_1, n, m)\\
     & \cdot \left[ \left(n\tau_0 - m \tau_1 - \lambda^-\tau_0\tau_1(n - m)\right)^2  + (\tau_0 + \tau_1 - 2\lambda^-\tau_0\tau_1)c(\tau_0, \tau_1, n, m) \right],
 \end{align*}
 where $C(\tau_0, \tau_1, n, m)$ and $c(\tau_0, \tau_1, n, m)$ are defined as follows
 \begin{align*}
    C(\tau_0, \tau_1, n, m) &\coloneqq \max\qty{(1 - \lambda^-\tau_0)^{-2n}(1-\lambda^-\tau_1)^{-2}, (1-\lambda^-\tau_0)^{-2}(1-\lambda^-\tau_1)^{-2m} }, \\
    c(\tau_0, \tau_1, n, m) &\coloneqq \min\left\{n(\tau_0 - \lambda^-\tau_0\tau_1), m(\tau_1 - \lambda^- \tau_0\tau_1) \right\}.
 \end{align*}   
\end{lemma}
\begin{proof}
    We first observe that 
    \begin{equation}
            \distance_Z^2(J_\tau^n z,z)\leq(n\tau)^2(1 - \lambda^-\tau)^{-2n}\slopephi^2(z),
            \label{eq:J_tau^n}
    \end{equation}
    for $z\in\Domain(\slopephi)$, $n\in\Z_{>0}$ and $\tau>0$.
    Indeed, \cref{lem:sloepe_estimate} yields that
    \begin{align*}
        \distance^2_Z(J_\tau^n z,z)&\le n\sum_{k=1}^n\distance_Z^2(J_\tau^{k}z,J_\tau^{k-1}z)\\
        &\le n\tau^2(1+\lambda\tau)^{-2}\sum_{k=1}^n\slopephi^2(J_\tau^{k-1}z)\\
        &\le n\tau^2(1+\lambda\tau)^{-2}\qty(\sum_{k=1}^n(1+\lambda\tau)^{-2(k-1)})\slopephi^2(z)\\
        &\le n^2\tau^2(1 - \lambda^-\tau)^{-2n}\slopephi^2(z).
    \end{align*}
    Next, for $z \in \Domain(\slopephi)$ and $\tau_0, \tau_1 \in \left(0, \frac{1}{\lambda^-} \right)$, we obtain that
    \begin{equation}
    \label{eq:J_tau0_J_tau1}
        \distance_Z^2(J_{\tau_0} z, J_{\tau_1}z)\le\frac{\tau_0 - \lambda^-\tau_0\tau_1}{\tau_0+\tau_1 - 2\lambda^-\tau_0\tau_1}\distance^2_Z(J_{\tau_0}z,z)+\frac{\tau_1- \lambda^-\tau_0\tau_1}{\tau_0+\tau_1 - 2\lambda^-\tau_0\tau_1}\distance^2_Z(z,J_{\tau_1}z).
    \end{equation}
    Indeed, using \cref{lem:discrete_VI} again, we have
    \begin{align*}
        &\frac{1}{2\tau_0}\qty(\distance_Z^2(J_{\tau_0}z,z'_0)-\distance_Z^2(z,z'_0))+\frac{\lambda}{2}\distance_Z^2(J_{\tau_0}z,z'_0)\le\phi(x'_0,J_{\tau_0}^Yz)-\phi(J_{\tau_{0}}^Xz,y'_0),\\
        &\frac{1}{2\tau_1}\qty(\distance_Z^2(J_{\tau_1}z,z'_1)-\distance_Z^2(z,z'_1))+\frac{\lambda}{2}\distance_Z^2(J_{\tau_1}z,z'_1)\le\phi(x'_1,J_{\tau_1}^Yz)-\phi(J_{\tau_{1}}^Xz,y'_1),
    \end{align*}
    for $z'_0$, $z'_1\in\Domain\phi$.
    Letting $z'_0=J_{\tau_1}z$ and $z'_1=J_{\tau_0}z$ and summing up the above gives \eqref{eq:J_tau0_J_tau1}.
    Finally, combining \eqref{eq:J_tau^n} and \eqref{eq:J_tau0_J_tau1} inductively with \cite[Lemma~A2]{clement2009introduction} yields the desired inequality.
\end{proof}

\begin{proposition}[Existence for the case $w_0 \in \Domain(\slopephi)$]\label{prop:convergence}
    Suppose that \cref{ass:function_value_on_and_out_domain,assump:countinuous_coupling_mild,assump:convex_concavity_Phi} hold for some $\lambda \in \R$.
    Then for any $z \in \Domain(\slopephi)$ and $t \ge 0$ there exists the limit $w_t\coloneqq\lim_{n\to\infty}J^n_{t/n}z$, where we define $J_0z \coloneqq z$.
    Moreover, the following hold:
    \begin{enumerate}[label=(\roman*),ref=(\roman*)]
        \item \label{prop:convergence:uniform} $J_{t/n}^nz$ converges to $w_t$ uniformly on $[0, T]$ for any $T > 0$;
        \item\label{prop:convergence:error} let $T > 0 $ and $N \geq 1$ such that $T/N \in \left(0, \frac{1}{\lambda^-}\right)$. Then the following error estimate holds
            \begin{align*}
                \distance_Z^2(w_t,J_{t/n}^nz)
                \le \slopephi^2(z)C(T, N, \lambda^-)\frac{1}{n} \text{ for any } t \in [0, T] \text{ and } n \geq N,
            \end{align*}
            where $C(T, N, \lambda^-) = (1 - \lambda^- T /N)^{-2(N+1)}T^2\left(1 + (T\lambda^-)^2/N\right)$;
        \item\label{prop:convergence:lipschitz} the map $[0, +\infty) \ni w_t$ is locally Lipschitz;
        \item\label{prop:convergence:unique} the map $(0, +\infty) \ni t \mapsto w_t$ is the unique solution of $\EVSI{X}{Y}{\phi}{\lambda}$ with $\lim_{t \downarrow 0}w_t=z$.
    \end{enumerate}
\end{proposition}
The proof is given by a similar argument to that of Crandall--Ligget and Clément in \cite{MR287357,clement2009introduction}.
In their work, the method was applied to nonlinear semigroups on Banach spaces and gradient flows on metric spaces, but as you can see from the following proof, it also works with our EVIs-flow on metric spaces.
\begin{proof}
Let $z\in\Domain(\slopephi)$, $t>0$ and $n, m \in \N$ be sufficiently large.
Letting $\tau_0 \coloneqq \nicefrac{T}{n}$ and $\tau_1 \coloneqq \nicefrac{t}{m}$ in \cref{lem:distance_between_resolvents} yields that
\begin{align*}
    \label{eq:Cauchy_estimate}
    \begin{aligned}
    \distance_Z^2(J_{t/m}^m z, J_{t/n}^n z)
    \le{}& \slopephi^2(z) C(t/n, t/m, n, m) \\
    & \cdot t^2 \qty(\qty(t\lambda^-\qty(\frac{1}{m} - \frac{1}{n}))^2 + \qty(\frac{1}{n} + \frac{1}{m} - 2t\lambda^-\frac{1}{n}\frac{1}{m})c(t/n, t/m, n, m)).
\end{aligned}
\end{align*}
Therefore, $(J_{t/n}^nz)_n\subset X\times Y$ is a locally uniformly Cauchy sequence.
Thus we can define $w\colon[0,+\infty)\longrightarrow X\times Y$ by $w_t\coloneqq\lim_{n\to\infty}J_{t/n}^nz$ if $t>0$ $w(t)\coloneqq z$ if $t=0$ with the error estimate
\begin{equation}
    \label{eq:convergence_estimate}
    \distance_Z^2(w_t,J_{t/n}^nz)
    \le \slopephi^2(z)\max\qty{\qty(1 - \lambda^-{t}/{n})^{-2n}, \qty(1 - \lambda^-t/n)^{-2}\e^{2\lambda^-t}} t^2 \qty(1 + (t\lambda^-)^2/n)\frac{1}{n},
\end{equation}
for sufficiently large $n \in \N$.
This completes the proof of the existence of $w_t$, \ref{prop:convergence:uniform} and \ref{prop:convergence:error}.

Next, we prove \ref{prop:convergence:lipschitz}.
Similarly, \cref{lem:distance_between_resolvents} induces the local Lipschitz continuity of $w$:
\begin{align*}
    \distance_Z(w_{t_1}, w_{t_2}) \le |t_1 - t_2|\e^{\lambda^- t_2}\slopephi(z)
    \text{ for any } t_1, t_2 \ge 0 \text{ with } t_1 < t_2.
\end{align*}
This completes the proof of \ref{prop:convergence:lipschitz}.

Finally, we show \ref{prop:convergence:unique}.
Let $(u_t, v_t) \coloneqq w_t$.
\cref{lem:sloepe_estimate} and the lower semicontinuity of $\slopephi$ induce that
\begin{align*}
    \slopephi(w_t) \leq \e^{-\lambda t} \slopephi(z) \text{ for any } t \ge 0;
\end{align*}
in particular $w_t \in \Domain\phi$ at every $t \ge 0$.
Thanks to \cref{assump:countinuous_coupling_mild} it is enough to show that 
\begin{align*}
    \frac{1}{2}\left(\distance_X^2(u_{t_2}, x) - \distance_X^2(u_{t_1}, x)  \right) &\le \int_{t_1}^{t_2} \phi(x, v_r) - \phi(u_r, v_r) \dd{r} - \frac{\lambda}{2} \int_{t_1}^{t_2} \distance_X^2(u_r, x)\dd{r}, \\
    \frac{1}{2}\left(\distance_Y^2(v_{t_2}, y) - \distance_Y^2(v_{t_1}, y)  \right) &\le \int_{t_1}^{t_2} \phi(u_r, v_r) - \phi(u_r, y) \dd{r} - \frac{\lambda}{2} \int_{t_1}^{t_2} \distance_Y^2(v_r, y)\dd{r}
\end{align*}
hold for any $t_1, t_2 \in (0, +\infty)$ with $t_1 < t_2$ and $(x, y) \in \Domain \phi$.
We prove only the first inequality; the other hand follows by a similar argument.
Let $T > t_2$ and $n \in \N$. We define the following three:
\begin{itemize}[wide=0pt]
    \item for any $t \in [0, T]$ we define $k_n(t) \in \{1, \cdots, n\}$ by $t \in ((k_n(t) - 1)\cdot\nicefrac{T}{n}, k_n(t)\cdot\nicefrac{T}{n}]$ if $t > 0$ and $k_n(t) = 0$ otherwise;
    \item we define $(\bar{x}_t^n, \bar{y}^n_t) \colon [0, T] \to \XY$ by $(\bar{x}_t^n, \bar{y}^n_t) \coloneqq J_{T/n}^{k_n(t)}z$;
    \item $I_m \coloneqq ((m-1)\cdot\nicefrac{T}{n}, m\cdot\nicefrac{T}{n}]$ for $m = 1, \ldots, n$.
\end{itemize}
For $(x^m, y^m) \coloneqq J_{T/n}^mz$ \cref{lem:discrete_VI} implies that
\begin{align*}
    \frac{1}{2}\distance_X^2(x^m, x) - \frac{1}{2}\distance_X^2(x^{m-1}, x)
    \le \frac{T}{n} \left( \phi(x, y^m) - \phi(x^m, y^m) - \frac{\lambda}{2}\distance_X^2(x^m, x) \right);
\end{align*}
thus we have
\[
    \begin{aligned}
    \frac{1}{2}\distance_X^2(\bar{x}^n_{t_2}, x) - \frac{1}{2}\distance_X^2(\bar{x}^n_{t_1}, x)
    &{}\le \sum_{m=k_n(t_1)+1}^{k_n(t_2)} \left( \int_{I_m} \phi(x, \bar{y}^n_r) - \phi(\bar{x}^n_r, \bar{y}^n_r) \dd{r} - \int_{I_m} \frac{\lambda}{2} \distance_X^2(\bar{x}^n_r, x) \dd{r} \right) \\
    &{}= \int_{k_n(t_1)T/n}^{k_n(t_2)T/n} \phi(x, \bar{y}^n_r) - \phi(\bar{x}^n_r, \bar{y}^n_r) \dd{r} - \frac{\lambda}{2} \int_{k_n(t_1)T/n}^{k_n(t_2)T/n} \distance_X^2(\bar{x}_r^n, x) \dd{r}.
    \end{aligned}
\]
By \cref{lem:distance_between_resolvents} and the similar argument of the proof of (i) it holds that $(\bar{x}^n_t, \bar{y}^n_t)$ uniformly converges to $w_t$ on $[0, T]$ as $n \to+\infty$; thus we have
\begin{align*}
    \frac{\lambda}{2} \int_{k_n(t_1)T/n}^{k_n(t_2)T/n} \distance_X^2(\bar{x}_r^n, x) \dd{r}
    \to \frac{\lambda}{2} \int_{t_1}^{t_2} \distance_X^2(u_r, x) \dd{r}.
\end{align*}
In addition we show that there exists $N_0 \in \N$ and $C > 0$ such that  $\sup_{n \ge N_0} \sup_{r \in [0, T]} \phi(x, \bar{y}^n_r) - \phi(\bar{x}^n_r, \bar{y}^n_r) \leq C$ in what follows; therefore by Fatou's lemma we can get the desired results i.e.~we completes the proof.
Suppose that $\limsup_{n \to+ \infty} \sup_{n \ge N_0} \sup_{r \in [0, T]} \phi(x, \bar{y}^n_r) - \phi(\bar{x}^n_r, \bar{y}^n_r) = +\infty$.
After passing to a subsequence we can take $(t_n)_{n\ge1} \in [0, T]$ such that $\limsup_{n\to +\infty} \phi(x, \bar{y}^n_{t_n}) - \phi(\bar{x}^n_{t_n}, \bar{y}^n_{t_n}) = + \infty$.
Since $[0, T]$ is compact, after passing to a subsequence, we have $t \in [0, T]$ such that $t_n \to t$ as $n \to+\infty$.
By \cref{lem:distance_between_resolvents} and the similar argument of the proof of (i) we have $\distance_Z\left((\bar{x}^n_{t_n}, \bar{y}^n_{t_n}), J_{t/n}^nz  \right)\to 0$ as $n \to+\infty$; thus $(\bar{x}^n_{t_n}, \bar{y}^n_{t_n}) \to w_t$ as $n \to+\infty$.
Since $(\bar{x}^n_{t_n}, \bar{y}^n_{t_n}) \in \Domain\phi$ for each $n \ge 1$ thanks to \cref{assump:countinuous_coupling_mild} we have
\begin{align*}
    +\infty &= \limsup_{n\to +\infty} \left( \phi(x, \bar{y}^n_{t_n}) - \phi(\bar{x}^n_{t_n}, \bar{y}^n_{t_n}) \right)\\
    &= \limsup_{n\to +\infty} \left( \varphi(x, \bar{y}^n_{t_n}) + \psi_X(x) - \varphi(\bar{x}^n_{t_n}, \bar{y}^n_{t_n}) - \psi_X(\bar{x}^n_{t_n}) \right) \\
    &\le \varphi(x, v_t) + \psi_X(x) - \varphi(u_t, v_t) - \psi_X(u_t).
\end{align*}
This implies that $\psi_X(u_t) = -\infty$; however this is a contradiction.
\end{proof}

We extend the above result to the case $w_0 \in \overline{\Domain \phi}$ under \cref{assump:countinuous_coupling}, which is stronger than \cref{assump:countinuous_coupling_mild}.
\begin{proof}[Proof of \cref{thm:existence}]
For each $t\ge0$ we define $S_t \colon \Domain(\slopephi) \to \Domain(\slopephi)$ to be that $S_t[z] \coloneqq w_t$ for any $z \in \Domain(\slopephi)$, where $w$ is the unique solution of $\EVSI{X}{Y}{\phi}{\lambda}$ with $w_0 = z$.
Then $\{S_t\}_{t \geq 0}$ is a $C_0$-semigroup thanks to \cref{thm:contraction}.
Moreover from \cref{thm:contraction} we can extend $\{S_t\}_{t\ge0}$ over $\overline{\Domain\phi}$ as a $C_0$-semigroup in an obvious way.
Let $(u_t, v_t)=w_t=S_t[z]$ for $z \in \overline{\Domain\phi}$ and $t\ge0$.
Since $\{S_t\}_{t\ge0}$ is $C_0$-semigroup (in particular $\lim_{t\downarrow 0} w_t = z$) and $\lambda$-contraction holds, the map $w_t \colon [0, +\infty) \to \overline{\Domain\phi}$ is continuous.
Therefore we can check that the maps $t \mapsto \phi(u_t, y), \phi(x, v_t), \phi(u_t, v_t)$ belong to $\Loneloc((0, +\infty))$ and \eqref{eq:EVIs'} hold for each $(x, y) \in \Domain\phi$, by taking an approximate sequence $\{z_n\}_{n\ge1} \subset \Domain(\slopephi)$ for $z$.
\cref{prop:integral_characterization} (b) completes the proof.
\end{proof}

\begin{remark}\label{rem:discussion_about_continuous_coupling}
    Under \cref{assump:countinuous_coupling_mild} rather than \cref{assump:countinuous_coupling}, we can derive the following results:
    for any  $z \in \overline{\Domain \phi}$ there exists a continuous map $w_t \colon [0, +\infty) \to \overline{\Domain \phi}$ satisfying \cref{prop:integral_characterization} (ii) and  $w_0 = z$. 

    Furthermore, $w_t$ is locally Lipschitz on $(0, +\infty)$ if and only if there exists a sequence $(t_n)_{n\ge1} \subset (0, 1)$ such that $t_n \to 0$ and $w_t$ is rectifiable around at each $t_n$.
    This follows from that $\delta_+(t_n) < +\infty$, \eqref{eq:regularlizing_effect_1} and the similar argument as the proof of \cref{prop:integral_characterization}.
    In particular, if $w_t$ is locally Lipschitz, then $w_t \in \Domain \phi$ for each $t > 0$ i.e. $w_t$ is the solution of $\EVSI{X}{Y}{\phi}{\lambda}$ with $\lim_{t\to0}w_t = z$.
\end{remark}
\section{Applications for the well-posedness of Wasserstein gradient descent-ascent flows}\label{sec:WGDA}
This chapter presents a key application of our theoretical framework developed in preceding sections: exponential convergence of Wasserstein gradient descent-ascent flows (Wasserstein GDA).
In this article, we give a unified geometric viewpoint---\cref{assump:convex_concavity_Phi}---for the existence and convergence of Wasserstein GDA.

We begin by recalling the fundamental definitions of optimal‐transport theory.
We follow the notation of \cite[Chapter~7]{AGS}.
Let $\X$, $\Y$ be (possibly distinct) separable Hilbert spaces.
We will denote by $\Ptwo(\X)$ and $\Ptwo(\Y)$ the set of Borel probability measures with finite second moments on $\X$ and $\Y$, respectively.
That is,  
\(
    \int_{X}\norm{x}_X^2\dd{\mu(x)} < +\infty\text{ for }\mu\in\Ptwo(\X).
\)
The ($2$-)Wasserstein distance $W_2$ between $\mu_1$, $\mu_2\in\Ptwo(\X)$ is defined as
\begin{equation}
    W_2(\mu^1,\mu^2)\coloneqq\qty(\inf_{\bmu}\int_{\X\times\X}\norm{x-x'}_\X^2\dd{{\bmu}(x,x')})^{1/2},\label{eq:def_W2}
\end{equation}
where $\bmu$ runs over the set of joint probability measures or couplings between $\mu_1$ and $\mu_2$.
We call the complete metric space $(\Ptwo(\X),W_2)$ \emph{the ($2$-)Wasserstein space over $\X$}.
Denoted by $\Gamma(\mu_1, \ldots, \mu_n)$ the all multiple plans with marginals $\mu_i \in \Ptwo(\X)$ and denoted by $\Gamma_{\textup{o}}(\mu, \nu) \subset \Gamma(\mu, \nu)$ the all optimal plans between $\mu \in \Ptwo(\X)$ and $\nu \in \Ptwo(\X)$.
In what follows, we assume that $\Ptwo(\X) \times \Ptwo(\Y)$ is equipped with $\ell^2$-distance.

\subsection{Preliminaries: absolutely continuous curves in the $2$-Wasserstein space \texorpdfstring{$(\Ptwo(\X), W_2)$}{W2} }\label{subsec:AC_curve_in_P2}
To define Wasserstein GDA, this section is dedicated to recalling a characterization of absolutely continuous curves in $\Ptwo(\X)$ developed in~\cite{AGS}.
Most of the results presented in this section are quoted from~\cite[Chapter~8]{AGS}.

We start by defining the tangent space to $\Ptwo(\X)$ at $\mu \in \Ptwo(\X)$.
Let $\qty{e_n}_{n=1}^\infty$ be an orthonormal basis of $\X$.
Denoted by $\cyl(\X)$ the collections of all functions $\zeta \coloneqq \eta \circ \pi$ with $\eta \in C_c^\infty(\R^n)$ and $\pi \colon \X \to \R^n$ of the form 
\begin{align*}
    \pi(x) \coloneqq \qty(\la x, e_1\ra, \ldots, \la x, e_n \ra ) \text{ for each } x \in \X,
\end{align*}
for each $n \in \N$.
Note that any $\zeta \in \cyl(\X)$ is Fr\'echet differentiable: we denote by $\nabla \zeta$ its Fr\'echet derivative, and the map $\X \ni x \mapsto \nabla \zeta(x)$ is continuous.
Using $\cyl(\X)$ as test functions, we define the tangent space to $\Ptwo(\X)$ at each $\mu \in \Ptwo(\X)$.
Let $L_2(\mu, \X)$ be the $L^2$-space which is consisted of $\X$-valued $\mu$-measurable functions.
We introduce the tangent space $\Tan_\mu \Ptwo(\X)$ to be that
\begin{align*}
    \Tan_\mu \Ptwo(\X) \coloneqq \overline{\{\nabla \zeta \mid \zeta \in \cyl(\X) \}}^{L^2(\mu, \X)}.
\end{align*}
For $v \in L^2(\mu, \X)$ let $S(v) \coloneqq \{ v + w \mid w \in L^2(\mu, \X) \text{ s.t. } \nabla \cdot w\mu = 0 \}$, where $\nabla \cdot w \mu = 0$ is in the sense of distribution with respect to $\cyl(\X)$:
\begin{align*}
    \int_{\X} \la w, \nabla \zeta \ra \dd{\mu} = 0 \text{ for any } \zeta \in \cyl(\X). 
\end{align*}
Since the norm $\| \cdot \|_{L^2(\mu, \X)}$ is strictly convex, it has the unique minimizer $m(v)$ in $S(v)$.
Then it holds that
\begin{align}\label{eq:equivalent_cond_tan_vector}
    v \in L^2(\mu, \X) \text{ belongs to } \Tan_\mu\Ptwo(\X) \text{ if and only if } v = m(v).
\end{align}

We now move to a general metric space.
Let $I \subset \R$ be an open interval and $X$ be a complete metric space.
A curve $u_t$ from $I$ to $X$ is said to be a \emph{absolutely continuous} (resp.~\emph{locally absolutely continuous}) if there exists $m \in L^1(I)$ (resp.~$\Loneloc(I)$) such that
\begin{align*}
    \distance_X(u_t, u_s) \leq \int_s^t m(r) \dd{r}\text{ for each } s, t \in I \text{ with } s < t.
\end{align*}
For a locally absolutely continuous curve $u_t$ it holds that $|u'_t| \coloneqq \lim_{s\to t}\frac{\distance_X(u_t, u_s)}{|t-s|}$ exists for $\mathscr{L}^1$-a.e.~$t \in I$.

We are now back in the Wasserstein space.
Locally absolutely continuous curves in $\Ptwo(\X)$  have the following characterization.
Let $\mu_t$ be a locally absolutely continuous curve from $I$ to $\Ptwo(\X)$.
Then there exists a Borel vector field $v \colon (x,t) \mapsto v_t(x)$ from $I \times \X$ to $\X$ satisfying 
\begin{align}\label{eq:tangent_vec_AC}
    v_t \in \Tan_{\mu_t} \Ptwo(\X) \text{ and }
    \|v_t\|_{L^2(\mu, \X)} = |\mu'_t|\text{ for } \mathscr{L}^1\text{-a.e. } t \in I,
\end{align}
and the continuity equation $\partial_t\mu_t + \nabla \cdot v_t\mu_t = 0$ which is defined to be that
\begin{align}\label{eq:cont_eq}
    \int_I\int_\X \partial_t\zeta(t, x) \dd \mu_t + \la v_t(x), \nabla_x \zeta(t, x)\ra \dd \mu_t = 0 \text{ for each } \zeta \in C_c^\infty(I) \times \cyl(\X).
\end{align}
Moreover, if a vector field $w$ also satisfies \eqref{eq:tangent_vec_AC} and \eqref{eq:cont_eq} for $\mu_t$, then it follows from \eqref{eq:equivalent_cond_tan_vector} that
\begin{align*}
    w_t = v_t \text{ for } \mathscr{L}^1\text{-a.e. } t \in I.
\end{align*}
In this sense, the vector field $v$ is uniquely determined for $\mathscr{L}^1$-a.e. $t \in I$.
We call $v_t$ the tangent vector for $\mu_t$ and denote $v_t$ by $\mu'_t$.

\subsection{Preliminaries: sub-differential calculus on \texorpdfstring{$(\Ptwo(\X), W_2)$}{W2} }\label{subsec:sub_differential_on_P2}
For the same purpose of \cref{subsec:AC_curve_in_P2} we recall the sub-differential calculus on $\Ptwo(\X)$.
Similarly to the previous section, most of the results presented in this section are quoted from~\cite[Chapter~10]{AGS}.
Let $f \colon \Ptwo(\X) \longrightarrow \overR$ be a proper, lower semi-continuous function and $\mu \in \Domain f$.
We say that $v \in L^2(\mu, \X)$ is a \textit{Fr\'echet sub-differential} of $f$ at $\mu$ if it holds that
\begin{align*}
    \liminf_{\Domain f \ni \nu \to \mu} \frac{f(\nu) - f(\mu) - \inf_{\bmu \in \Gamma_{\textup{o}}(\mu, \nu)} \int_{\X^2} \la v(x), y - x\ra \dd{\bmu(x,y)}} {W_2(\mu, \nu)} \geq 0.
\end{align*}
Denoted by $\fsubd f(\mu)$ the all Fr\'echet sub-differentials of $f$ at $\mu$.
In the next section, we consider curves maximizing a concave function.
Unfortunately, this sub-differential notion is unsuitable for this purpose: for example, the concave function $\phi(x) \coloneqq - |x|$ on $\mathbb{R}$ is not Fr\'echet sub-differentiable at $x=0$.
Thus, we also consider another differential notion called \textit{Fr\'echet super-differential}.
Let $g \colon \Ptwo(\X) \longrightarrow [-\infty, +\infty)$ be proper upper semi-continuous function and $\tilde{\mu} \in \Domain g  \coloneqq \{ g > -\infty\}$.
We say that $v \in L^2(\tilde{\mu}, \X)$ is a \textit{Fr\'echet superdifferential} of $f$ at $\tilde{\mu}$ if it holds that
\begin{align*}
    \limsup_{\Domain g \ni \tilde{\nu} \to \tilde{\mu}}\frac{g(\tilde{\nu}) - g(\tilde{\mu}) -  \sup_{\tilde{\bmu} \in \Gamma_{\textup{o}}(\tilde{\mu}, \tilde{\nu})} \int_{\X^2} \la v(x), y - x \ra \dd{\tilde{\bmu}(x, y)}}{W_2(\tilde{\mu}, \tilde{\nu})} \leq 0.
\end{align*}
Denoted by $\fsupd g(\tilde{\mu})$ the all Fr\'echet super-differentials of $g$ at $\tilde{\mu}$.
By the definition we can easily check that $v \in \fsubd f(\mu)$ if and only if $-v \in \fsupd(-f)(\mu)$.

Let $\lambda \in \R$ and  $f$ be also $\lambda$-convex along geodesics, which is defined as that, for any $\mu_0, \mu_1 \in \Ptwo(\X)$, there exists a geodesic $\mu_t \colon [0, 1] \to \Ptwo(\X)$ from $\mu_0$ to $\mu_1$ such that
\begin{align*}
    f(\mu_t) \leq (1-t) f(\mu_0) + t f(\mu_1) - \frac{\lambda}{2} W_2^2(\mu_0, \mu_1) \text{ for each } t \in [0, 1].
\end{align*}
Similarly, as for $\lambda$-convex functionals on Banach spaces, the following holds: 
a vector $v \in L^2(\mu, \X)$ belongs to $\fsubd f(\mu)$ if and only if for any $\nu \in \Ptwo(\X)$ there exists $\gamma \in \Gamma_{\textup{o}}(\mu, \nu)$ such that
\begin{align}\label{eq:equivalence_subdifferential_inclusion}
    f(\nu) - f(\mu) \geq \int_{\X^2} \la v(x), y - x \ra \dd{\gamma(x, y)}+ \frac{\lambda}{2}W_2^2(\mu, \nu).
\end{align}

\subsection{Existence and convergence of Wasserstein GDA}\label{subsec:conv_WGDA}

Let $\phi \colon \Ptwo(\X) \times \Ptwo(\Y) \longrightarrow \overline{\R}$ be a closed functional.
We define Wasserstein GDA as follows:
\begin{definition}[Wasserstein gradient descent-ascent flows]
    \label{def:WassGDA}
    A pair of locally absolutely continuous curves $(\mu,\nu)\colon \zeroinfty \longrightarrow\Ptwo(\X)\times\Ptwo(\Y)$ is said to be a \emph{Wasserstein gradient descent-ascent flow} if it satisfy 
    \[
    \left\{
    \begin{aligned}
        -\mu'_t & \in \fsubd_\X \phi(\mu_t, \nu_t), \\
        \nu'_t & \in \fsupd_\Y \phi(\mu_t, \nu_t), 
    \end{aligned}
    \right.
    \]
    for $\Lebae t \in \zeroinfty$,
    where $\mu_t'$ (resp.~$\nu_t'$) is the tangent vector for $\mu_t$ (resp.~$\nu_t$), and $\fsubd_\X$ (resp.~$\fsupd_\Y$) is the Fr\'echet sub-differential (resp.~super-differential) for the functional $\Ptwo(\X) \ni \mu \mapsto \phi(\mu, \nu_t)$ (resp.~$\Ptwo(\Y) \ni \nu \mapsto \phi(\mu_t, \nu)$).
\end{definition}
To establish a link between EVIs and Wasserstein GDA, we introduce a convexity-concavity along a specific curve called \textit{generalized geodesic}.
First, we recall the generalized geodesics introduced in~\cite[Definition~9.2.2]{AGS}.
Let $\mu_b$, $\mu_0$, $\mu_1 \in \Ptwo(\X)$ and $\bmu \in \Gamma(\mu_b, \mu_0, \mu_1)$ be satisfying $\pi^{1, 2}_\#\bmu \in \Gamma_{\textup{o}}(\mu_b, \mu_0)$ and $\pi^{1, 3}_\#\bmu \in \Gamma_{\textup{o}}(\mu_b, \mu_1)$, where $\pi^i \colon \prod_{i=1}^3 \X \to \X$ is the canonical projection.
A generalized geodesic $\bmu_t \colon [0, 1] \to \Ptwo(\X)$ from $\mu_0$ to $\mu_1$ with base $\mu_b$ is defined to be that
\begin{align*}
    \bmu_t \coloneqq (\pi^{2\to3}_t)_\#\bmu,
\end{align*}
where $\pi^{2\to3}_t \coloneqq (1-t) \pi^2 + t \pi^3$.
We also define
\begin{align*}
    W_\bmu(\mu_0, \mu_1)^2 \coloneqq \int_{\X^3} \|x_2 - x_3\|^2 \dd \bmu(x_1, x_2, x_3).
\end{align*}
Next, we introduce a convexity-concavity along generalized geodesics.
\begin{definition}[Modified $\lambda$-convex-concave along generalized geodesics]\label{def:modified_lambda_convex_concave}
Let $\mu_0, \mu_1, \mu_b \in \Ptwo(\X)$ and $\nu_0, \nu_1, \nu_b \in \Ptwo(\Y)$, and let $\bmu_t$ and $\bnu_t$ be generalized geodesics with base points $\mu_b$ and $\nu_b$, respectively.
The functional $\phi$ is said to be \emph{modified $\lambda$-convex-concave} along $\bmu_t$ and $\bnu_t$ if the following two conditions are satisfied:
\begin{enumerate}
    \item for any $\nu \in \Ptwo(\Y)$ with $\max\{\phi(\mu_0, \nu), \phi(\mu_1, \nu) \} < +\infty$, the functional $\Ptwo(\X) \ni \mu \mapsto \phi(\mu, \nu)$ is \textit{modified $\lambda$-convex} along with $\bmu_t$ , which is defined by
        \begin{align}\label{eq:modified_lambda_convex}
            \phi(\mu_t, \nu) \leq (1-t) \phi(\mu_0, \nu) + t \phi(\mu_1, \nu) - \frac{\lambda}{2} W_{\bmu}^2(\mu_0, \mu_1) \text{ for each } t \in [0, 1].
        \end{align}
    \item for any $\mu \in \Ptwo(\X)$ with $\min\{\phi(\mu, \nu_0), \phi(\mu, \nu_1) \} > -\infty$, the map $\Ptwo(\Y) \ni \nu \mapsto -\phi(\mu, \nu)$ is modified $\lambda$-convex along $\bnu_t$.
\end{enumerate}
\end{definition}
We also say that $\phi$ is modified $\lambda$-convex-concave along with generalized geodesics if the following holds:
for any $\mu_0, \mu_1, \mu_b \in \Ptwo(\X)$ and $\nu_0, \nu_1, \nu_b \in \Ptwo(\Y)$, there exist generalized geodesics $\bmu_t$ and $\bnu_t$ with base points $\mu_b$ and $\nu_b$ respectively such that $\phi$ is modified $\lambda$-convex-concave along with $\bmu_t$ and $\bnu_t$.
\begin{remark}
    In~\cite[Definition~9.2.4]{AGS} the functional $\Ptwo(\X) \ni \mu \mapsto f(\mu)$ satisfying \eqref{eq:modified_lambda_convex} is called as \textit{$\lambda$-convex} simply; however, in this paper, we distinguish it from \cref{def:lambda_convex_concave_along_curve} and adopt the term ``modified'' convexity along generalized geodesics.
\end{remark}
\begin{remark}\label{rem:modified_convex_concave_implies_convex_concave}
    When $\lambda \neq 0$, 
    \cref{def:modified_lambda_convex_concave} is slightly different from \textit{$\lambda$-convex-concave along generalized geodesics} in the sense of~\cref{def:lambda_convex_concave_along_curve}, because it holds that
    \begin{align*}
        W_\bmu(\mu_1, \mu_2) \geq W_2(\mu_1, \mu_2).
    \end{align*}
    In general, the equality above does not hold.
    From the above inequality, the following hold:
    \begin{enumerate}
        \item if $\lambda \geq 0$ and $\phi$ is modified $\lambda$-convex-concave along with generalized geodesics, then $\phi$ is also $\lambda$-convex-concave along with generalized geodesics;
        
        \item if $\lambda \leq 0$ and $\phi$ is $\lambda$-convex-concave along with generalized geodesics, then $\phi$ is also modified $\lambda$-convex-concave along with generalized geodesics.
    \end{enumerate}
\end{remark}
\begin{remark}
    In particular, letting $\mu_b \coloneqq \mu_0$, $\bmu_t$ is the geodesic from $\mu_0$ to $\mu_1$ and it holds that
    \begin{align*}
        W_\bmu(\mu_0, \mu_1) = W_2(\mu_0, \mu_1).
    \end{align*}
    Therefore, if $\phi$ is modified $\lambda$-convex-concave along with generalized geodesics, then $\phi$ is $\lambda$-convex-concave along with geodesics in the sense of~\cref{def:lambda_convex_concave_along_curve}.
\end{remark}
Modified $\lambda$-convex-concavity plays a key role in this article.
For any $\mu \in \Ptwo(\X)$ the map $\Ptwo(\X) \ni \mu' \mapsto  W_2^2(\mu', \mu)$ is modified $1$-convex along with generalized geodesics with the base point $\mu$; thus we get the next result immediately.
\begin{remark}
    If $\phi$ is modified $\lambda$-convex-concave along with generalized geodesics for some $\lambda \in \R$, then
    for any $(\mu, \nu) \in \Ptwo(\X) \times \Ptwo(\Y)$ the functional $\Ptwo(\X) \times \Ptwo(\Y) \ni (\mu', \nu') \mapsto \Phi_\tau(\mu, \nu; \mu', \nu')$ is modified $(\tau^{-1} + \lambda)$-convex-concave along with generalized geodesics for any $\tau \in \qty(0, \frac{1}{\lambda^-})$.
    In particular, from \cref{rem:modified_convex_concave_implies_convex_concave} we can easily check that $\phi$, $\Ptwo(\X)$ and $\Ptwo(\Y)$ satisfy \cref{assump:convex_concavity_Phi} for $\lambda$.
\end{remark}

We are now able to establish the equivalence of EVIs and Wasserstein GDA.
The proof of the following theorem is similar to \cite[Theorem~11.1.4]{AGS}.
\begin{theorem}
    \label{thm:equivalence}
    Suppose that \cref{assump:countinuous_coupling} holds and $\phi$ is also modified $\lambda$-convex-concave along with generalized geodesics for some $\lambda \in \R$.
    Let $(\mu_t, \nu_t) \colon [0, +\infty)$ be a curve with $(\mu_0, \nu_0) \in \overline{\Domain \phi}$.
    Then, the following are equivalent:
    \begin{enumerate}
        \item the map $(0, +\infty) \ni t \mapsto (\mu_t, \nu_t)$ is the solution of $\EVSI{\Ptwo(\X)}{\Ptwo(\Y)}{\phi}{\lambda}$;
        
        \item the map $(0, +\infty) \ni t \mapsto (\mu_t, \nu_t)$ is locally absolutely continuous and is a Wasserstein gradient descent-ascent flow.
    \end{enumerate}
\end{theorem}
\begin{proof}
    First, we show the direction: (i) $\Rightarrow$ (ii).
    By~\cref{thm:contraction} the maps $(0, +\infty) \ni t \mapsto (\mu_t, \nu_t)$ is locally absolutely continuous; thus there exist tangent vector $\mu'_t$ and $\nu_t'$ for $\Lebae t \in (0, +\infty)$.
    Let $\varphi$, $\psi_X$, and $\psi_Y$ be as in~\cref{assump:countinuous_coupling}.
    By~\cite[Theorem~8.4.7]{AGS} and \cref{def:EVIs}, for any countable sets $D_X \subset \Domain \psi_X $ and $D_Y \subset \Domain \psi_Y $, there exists a Lebesgue negligible set $\mathscr{N} \subset \zeroinfty$ such that,for any $t \in \zeroinfty \setminus \mathscr{N}$
    \begin{subequations}
    \begin{align}
        \label{eq:subdifferential_on_countable_element}
        &\int_{\X^2} \la -\mu_t'(x_1), x_2 - x_1 \ra \dd\bmu(x_1, x_2)  + \frac{\lambda}{2}\distance_X^2(\mu_t, \mu) \leq \phi(\mu, \nu_t) - \phi(\mu_t, \nu_t),\\
        \label{eq:superdifferential_on_countable_element}
        &\int_{\Y^2} \la -\nu_t'(y_1), y_2 - y_1 \ra \dd\bnu(y_1, y_2) + \frac{\lambda}{2}\distance_Y^2(\nu_t, \nu) \leq \phi(\mu_t, \nu_t) - \phi(\mu_t, \nu),
    \end{align}
    \label{eq:differential_ineqs}
    \end{subequations}
    for any $\mu \in D_X$, $\nu \in D_Y$, $\bmu \in \Gamma_{\textup{o}}(\mu_t, \mu)$ and $\bnu \in \Gamma_{\textup{o}}(\nu_t, \nu)$.
    In particular, we can take a Lebesgue negligible set $\mathscr{N}$ such above for countable dense subsets $D_X \subset \Domain\phiX$ and $D_Y \subset \Domain\phiY$ with respect to the distance $W_2(\mu^1, \mu^2) + |\phiX(\mu^1) - \phiX(\mu^2)|$ for $\mu^1$, $\mu^2\in D_X$ and $W_2(\nu^1, \nu^2) + |\phiY(\nu^1) - \phiY(\nu^2)|$ for $\nu^1$, $\nu^2\in D_Y$.
    For any $\mu \in \Domain\phiX$ (resp.~$\nu \in \Domain\phiY$) we can take a convergent sequence $(\mu_n)_{n} \subset D_X$ (resp.~$(\nu_n)_{n} \subset D_Y$); thus from \eqref{eq:subdifferential_on_countable_element} and \eqref{eq:superdifferential_on_countable_element}, noting that the continuity of $\contiphi$ and narrowly lower semicontinuity of first terms in~\eqref{eq:subdifferential_on_countable_element} and~\eqref{eq:superdifferential_on_countable_element}, one sees that the same inequalities as \eqref{eq:differential_ineqs} holds   for any $\mu \in \Domain \phiX $ and $\nu \in \Domain \phiY $.
    Since $\Domain\phi = \Domain \phiX  \times \Domain \phiY $, from~\eqref{eq:equivalence_subdifferential_inclusion}, $(\mu_t, \nu_t)$ is a Wasserstein gradient descent-ascent flow.

    Converse direction directly follows from~\eqref{eq:equivalence_subdifferential_inclusion} and \cite[Theorem~8.4.7]{AGS}.
\end{proof}

The following result is a direct consequence of \cref{thm:equivalence,thm:contraction,thm:existence}.
\begin{corollary}\label{cor:existence_of_wgda}
    Suppose that \cref{assump:countinuous_coupling} holds and $\phi$ is also modified $\lambda$-convex-concave along with generalized geodesics for some $\lambda > 0$.
    Then for any $(\mu_0, \nu_0) \in \overline{\Domain \phi}$ there exists a unique Wasserstein gradient descent-ascent flow $(\mu_t, \nu_t) \colon [0, +\infty) \to \Ptwo(\X) \times \Ptwo(\Y)$.
    Moreover, for the unique saddle point $(\mu^\ast, \nu^\ast) \in \Domain \phi$ it holds that
    \begin{align*}
        W_2^2(\mu_t, \mu^\ast) + W_2^2(\nu_t, \nu^\ast) \leq \e^{-2\lambda t} \qty(W_2^2(\mu_0, \mu^\ast) + W_2^2(\nu_0, \nu^\ast)) \text{ for any } t \geq 0.
    \end{align*}
\end{corollary}

For example, the assumptions of \cref{cor:existence_of_wgda} is satisfied by the following functional:
\begin{proposition}\label{prop:application_functional_is_local_lip}
    Let $\ell \colon \X \times \Y \to \R$ be $C^{1,1}$ with $L\coloneqq \Lip(D\ell)$, where $D\ell$ is the Fr\'echet derivative of $\ell$.
    Then, the map $\Ptwo(\X) \times \Ptwo(\Y) \ni (\mu, \nu) \mapsto \int \ell(x, y) \dd{(\mu\otimes\nu)(x,y)}$ is locally Lipschitz; in particular, it belongs to $\R$.
\end{proposition}
\begin{proof}
    Let $B_R \subset \Ptwo(\X)\times\Ptwo(\Y)$ be the ball with the radius $R > 0$ at the center $(\mu_\ast, \nu_\ast) \in \Ptwo(\X)\times\Ptwo(\Y)$.
    Let $\mu, \mu' \in \Ptwo(\X)$, $\nu, \nu' \in \Ptwo(\Y)$ and $\varPi$ be a coupling between $\mu \otimes \nu$ and $\mu' \otimes \nu'$.
    By simple calculations, we have
    \begin{equation}
        \begin{aligned}
        &\abs{\int \ell(x, y) \dd{(\mu\otimes\nu)(x,y)} - \int \ell(x, y) \dd{(\mu\otimes \nu)(x,y)}}\\
        \leq &\int \abs{\ell(x, y) - \ell(x', y')} \dd{\varPi(x, y, x', y')} \\
        \le & \int \abs{\int_0^1 \la D\ell((1-t)x + tx', (1-t)y + ty'), (x' -x, y' -y) \ra \dd{t} } \dd{\varPi(x, y, x', y')} \\
        \leq & \int \qty[ \norm{D\ell(x,y)} + \frac{\sqrt{2}}{2}L \left(\|x - x'\| + \|y - y'\|\right) ] \sqrt{\|x - x'\|^2 + \|y - y'\|^2} \dd{\varPi(x,y, x', y')},
    \end{aligned}
    \label{eq:loc_lip_1}
    \end{equation}
    where we use the triangle inequality for $D\ell$ in the third inequality.
    In addition, for some $(x_0, y_0) \in \X \times \Y$,  by the triangle inequality, we get
    \begin{align}
        \label{eq:loc_lip_2}
        \|D\ell(x, y) \|^2 \leq 2 \left[ 4L^2 \left( \|x - x_0\|^2 + \|y - y_0\|^2\right) + \|D\ell(x_0, y_0)\|^2 \right].
    \end{align}
    Letting $\bmu_{\textup{o}}$ (resp.~$\bnu_{\textup{o}}$) be an optimal coupling between $\mu$ and $\mu'$ (resp.~$\nu$ and $\nu'$), note that $\bmu_{\textup{o}} \otimes \bnu_{\textup{o}}$ is a coupling of $\mu \otimes \nu$ and $\mu' \otimes \nu'$ (under a suitable push-forward).
    Thus, combining \eqref{eq:loc_lip_1} and \eqref{eq:loc_lip_2}, we have
    \begin{align*}
        \abs{\int \ell(x, y) \dd{(\mu\otimes\nu)(x,y)} - \int \ell(x, y) \dd{(\mu'\otimes \nu')(x,y)}}
        \leq C \sqrt{W_2^2(\mu, \mu') + W_2^2(\nu, \nu')},
    \end{align*}
    where $C$ is defined by
    \begin{align*}
        C^2 = 2 \|D\ell(x_0, y_0)\|^2 + 32L^2 \left(2R^2 + \int \|x - x_0\|^2 \dd{\mu_\ast(x)} + \int \|y - y_0\|^2 \dd{\nu_\ast(y)} \right),
    \end{align*}
    which depends on $L$, \(D\ell(x_0, y_0)\), \(\mu_\ast\), \(\nu_\ast\) and \(R\).
\end{proof}
\makeatletter
\begingroup
  \let\addcontentsline\@gobblethree
  \section*{Declaration of generative AI and AI-assisted technologies in the writing process}
  During the preparation of this work the authors used ChatGPT and Grammarly in order  to correct typographical errors.
  After using these services, the authors reviewed and edited the content as needed and take full responsibility for the content of the publication.
\endgroup
\makeatother

\printbibliography

@book{AGS,
      author={Ambrosio, Luigi and Gigli, Nicola and Savar{\'e}, Giuseppe},
       title={Gradient flows in metric spaces and in the space of probability
  measures},
     edition={2},
      series={Lectures in Mathematics ETH Z\"urich},
   publisher={Birkh\"auser Verlag},
     address={Basel},
        year={2008},
}

@InProceedings{pmlr-v247-wang24c,
  title = 	 {Open problem: Convergence of single-timescale mean-field Langevin descent-ascent for two-player zero-sum games},
  author =       {Wang, Guillaume and Chizat, L{\'e}na\"{i}c},
  booktitle = 	 {Proceedings of Thirty Seventh Conference on Learning Theory},
  pages = 	 {5345--5350},
  year = 	 {2024},
  volume = 	 {247},
  series = 	 {Proceedings of Machine Learning Research},
  url = 	 {https://proceedings.mlr.press/v247/wang24c.html}
}

@article{Rockafellar76,
author = {Rockafellar, R. Tyrrell},
title = {Monotone Operators and the Proximal Point Algorithm},
journal = {SIAM Journal on Control and Optimization},
volume = {14},
number = {5},
pages = {877--898},
year = {1976},
doi = {10.1137/0314056},
URL = { 
        https://doi.org/10.1137/0314056
},
eprint = { 
        https://doi.org/10.1137/0314056
}
}

@inproceedings{NEURIPS2020_e97c864e,
 author = {Domingo-Enrich, Carles and Jelassi, Samy and Mensch, Arthur and Rotskoff, Grant and Bruna, Joan},
 booktitle = {Advances in Neural Information Processing Systems},
 pages = {20215--20226},
 publisher = {Curran Associates, Inc.},
 title = {A mean-field analysis of two-player zero-sum games},
 url = {https://proceedings.neurips.cc/paper_files/paper/2020/file/e97c864e8ac67f7aed5ce53ec28638f5-Paper.pdf},
 volume = {33},
 year = {2020}
}

@inproceedings{
kim2024symmetric,
title={Symmetric Mean-field Langevin Dynamics for Distributional Minimax Problems},
author={Juno Kim and Kakei Yamamoto and Kazusato Oko and Zhuoran Yang and Taiji Suzuki},
booktitle={The Twelfth International Conference on Learning Representations},
year={2024},
url={https://openreview.net/forum?id=YItWKZci78}
}

@article{Conforti23,
author = {Conforti, Giovanni and Kazeykina, Anna and Ren, Zhenjie},
title = {Game on Random Environment, Mean-Field Langevin System, and Neural Networks},
journal = {Mathematics of Operations Research},
volume = {48},
number = {1},
pages = {78--99},
year = {2023},
doi = {10.1287/moor.2022.1252},

URL = { 
    
        https://doi.org/10.1287/moor.2022.1252
    
    

},
eprint = { 
    
        https://doi.org/10.1287/moor.2022.1252
    
    

}
}

@article{MURATORI2020108347,
title = {Gradient flows and Evolution Variational Inequalities in metric spaces. I: Structural properties},
journal = {Journal of Functional Analysis},
volume = {278},
number = {4},
pages = {108347},
year = {2020},
issn = {0022-1236},
doi = {https://doi.org/10.1016/j.jfa.2019.108347},
url = {https://www.sciencedirect.com/science/article/pii/S0022123619303416},
author = {Matteo Muratori and Giuseppe Savaré},
}

@article{Peiyuan23,
author = {Zhang, Peiyuan and Zhang, Jingzhao and Sra, Suvrit},
title = {Sion’s Minimax Theorem in Geodesic Metric Spaces and a Riemannian Extragradient Algorithm},
journal = {SIAM Journal on Optimization},
volume = {33},
number = {4},
pages = {2885--2908},
year = {2023},
doi = {10.1137/22M1505475},
URL = { 
        https://doi.org/10.1137/22M1505475
},
eprint = { 
        https://doi.org/10.1137/22M1505475
}
}

@incollection {MR285942,
    AUTHOR = {Rockafellar, R. T.},
     TITLE = {Monotone operators associated with saddle-functions and
              minimax problems},
 BOOKTITLE = {Nonlinear {F}unctional {A}nalysis},
    SERIES = {Proc. Sympos. Pure Math.},
    VOLUME = {XVIII, Part 1},
     PAGES = {241--250},
 PUBLISHER = {Amer. Math. Soc., Providence, RI},
      YEAR = {1970},
   MRCLASS = {47.80},
  MRNUMBER = {285942},
MRREVIEWER = {E.\ Asplund},
}

@incollection {Rockafellar71,
    AUTHOR = {Rockafellar, R. Tyrrell},
     TITLE = {Saddle-points and convex analysis},
 BOOKTITLE = {Differential {G}ames and {R}elated {T}opics ({P}roc.
              {I}nternat. {S}ummer {S}chool, {V}arenna, 1970)},
     PAGES = {109--127},
 PUBLISHER = {North-Holland Publishing Co., Amsterdam-London, American
              Elsevier Publishing Co., Inc., New York},
      YEAR = {1971},
   MRCLASS = {47.90 (46.00)},
  MRNUMBER = {285947},
MRREVIEWER = {Ky\ Fan},
}

@Inbook{Arrow2014,
author="Arrow, Kenneth J.
and Hurwicz, Leonid",
title="A Gradient Method For Approximating Saddle Points and Constrained Maxima",
bookTitle="Traces and Emergence of Nonlinear Programming",
year="2014",
publisher="Springer Basel",
address="Basel",
pages="45--60",
isbn="978-3-0348-0439-4",
doi="10.1007/978-3-0348-0439-4_2",
url="https://doi.org/10.1007/978-3-0348-0439-4_2"
}

@article{doi:10.1080/02331934.2023.2215799,
author = {Simon K. Niederländer},
title = {On the Arrow--Hurwicz differential system for linearly constrained convex minimization},
journal = {Optimization},
volume = {73},
number = {7},
pages = {2313--2345},
year = {2024},
publisher = {Taylor \& Francis},
doi = {10.1080/02331934.2023.2215799},


URL = { 
    
        https://doi.org/10.1080/02331934.2023.2215799
    
    

},
eprint = { 
    
        https://doi.org/10.1080/02331934.2023.2215799
    
    

}

}

@InProceedings{pmlr-v238-dvurechensky24a,
  title = 	 {Analysis of Kernel Mirror Prox for Measure Optimization},
  author =       {Dvurechensky, Pavel and Zhu, Jia-Jie},
  booktitle = 	 {Proceedings of The 27th International Conference on Artificial Intelligence and Statistics},
  pages = 	 {2350--2358},
  year = 	 {2024},
  volume = 	 {238},
  series = 	 {Proceedings of Machine Learning Research},
  url = 	 {https://proceedings.mlr.press/v238/dvurechensky24a.html},
}

@article{CHERUKURI201610,
title = {Asymptotic convergence of constrained primal-dual dynamics},
journal = {Systems \& Control Letters},
volume = {87},
pages = {10--15},
year = {2016},
issn = {0167-6911},
doi = {https://doi.org/10.1016/j.sysconle.2015.10.006},
url = {https://www.sciencedirect.com/science/article/pii/S0167691115002078},
author = {Ashish Cherukuri and Enrique Mallada and Jorge Cortés}
}

@INPROCEEDINGS{9483346,
  author={You, Pengcheng and Mallada, Enrique},
  booktitle={2021 American Control Conference (ACC)}, 
  title={Saddle Flow Dynamics: Observable Certificates and Separable Regularization}, 
  year={2021},
  volume={},
  number={},
  pages={4817--4823},
  doi={10.23919/ACC50511.2021.9483346}}

@book{Ekeland99,
author = {Ekeland, Ivar and Témam, Roger},
title = {Convex Analysis and Variational Problems},
publisher = {Society for Industrial and Applied Mathematics},
year = {1999},
}

@Article{Kopfer2017,
author={Kopfer, Eva},
title={Gradient flow for the Boltzmann entropy and Cheeger's energy on time-dependent metric measure spaces},
journal={Calculus of Variations and Partial Differential Equations},
year={2017},
day={23},
volume={57},
number={1},
pages={20},
issn={1432-0835},
doi={10.1007/s00526-017-1287-5},
url={https://doi.org/10.1007/s00526-017-1287-5}
}

@Article{Ferreira2018,
author={Ferreira, Lucas C. F.
and Valencia-Guevara, Julio C.},
title={Gradient flows of time-dependent functionals in metric spaces and applications to {PDE}s},
journal={Monatshefte f{\"u}r Mathematik},
year={2018},
day={01},
volume={185},
number={2},
pages={231--268},
issn={1436-5081},
doi={10.1007/s00605-017-1037-y},
url={https://doi.org/10.1007/s00605-017-1037-y}
}

@article {Asakawa86,
    AUTHOR = {Asakawa, Hidekazu},
     TITLE = {Maximal monotone operators associated with saddle functions},
   JOURNAL = {TRU Math.},
  FJOURNAL = {TRU Mathematics},
    VOLUME = {22},
      YEAR = {1986},
    NUMBER = {2},
     PAGES = {47--71},
      ISSN = {0496-6597},
   MRCLASS = {90C25 (46N05 47H05 90C48)},
  MRNUMBER = {908064},
MRREVIEWER = {Jean-Pierre\ Crouzeix},
}

@article{MIMURA20171477,
title = {The variational formulation of the fully parabolic {Keller--Segel} system with degenerate diffusion},
journal = {Journal of Differential Equations},
volume = {263},
number = {2},
pages = {1477--1521},
year = {2017},
issn = {0022-0396},
doi = {https://doi.org/10.1016/j.jde.2017.03.020},
url = {https://www.sciencedirect.com/science/article/pii/S0022039617301559},
author = {Yoshifumi Mimura},
}

@article{Blanchet+15,
	author = {{Blanchet, Adrien} and {Carrillo, José Antonio} and {Kinderlehrer, David} and {Kowalczyk, Michał} and {Laurençot, Philippe} and {Lisini, Stefano}},
	title = {A hybrid variational principle for the Keller–Segel system in
            $\mathbb{R}^2$},
	DOI= "10.1051/m2an/2015021",
	url= "https://doi.org/10.1051/m2an/2015021",
	journal = {ESAIM: M2AN},
	year = 2015,
	volume = 49,
	number = 6,
	pages = "1553--1576",
	month = "",
}

@misc{conger2024coupledwassersteingradientflows,
      title={Coupled Wasserstein Gradient Flows for Min-Max and Cooperative Games}, 
      author={Lauren Conger and Franca Hoffmann and Eric Mazumdar and Lillian J. Ratliff},
      year={2024},
      eprint={2411.07403v1},
      archivePrefix={arXiv},
      primaryClass={math.AP},
      url={https://arxiv.org/abs/2411.07403}, 
}

@inproceedings{Conger23,
 author = {Conger, Lauren and Hoffmann, Franca and Mazumdar, Eric and Ratliff, Lillian},
 booktitle = {Advances in Neural Information Processing Systems},
 pages = {45971--46006},
 title = {Strategic Distribution Shift of Interacting Agents via Coupled Gradient Flows},
 url = {https://proceedings.neurips.cc/paper_files/paper/2023/file/902c462e821e5e639ac3422b48b65932-Paper-Conference.pdf},
 volume = {36},
 year = {2023}
}

@misc{clement2010introduction,
  title={Introduction to Gradient Flows in Metric Spaces (II)},
  author={Cl{\'e}ment, Philippe},
  year={2010},
  url={https://igk.math.uni-bielefeld.de/sites/default/files/study_materials/notes2.pdf},
}

@misc{clement2009introduction,
  title={Introduction to Gradient Flows in Metric Spaces},
  author={Cl{\'e}ment, Philippe},
  year={2009},
  url={https://math.leidenuniv.nl/reports/files/2009-09.pdf},
}

@article{GOSSEZ1972220,
title = {On the subdifferential of a saddle function},
journal = {Journal of Functional Analysis},
volume = {11},
number = {2},
pages = {220--230},
year = {1972},
issn = {0022-1236},
doi = {https://doi.org/10.1016/0022-1236(72)90092-4},
url = {https://www.sciencedirect.com/science/article/pii/0022123672900924},
author = {Jean-Pierre Gossez},
}

@article{Awi24,
author = {Awi, Romeo and Hynd, Ryan and Mawi, Henok},
title = {Continuous time approximation of Nash equilibria in monotone games},
journal = {Communications in Contemporary Mathematics},
volume = {26},
number = {05},
pages = {2350021},
year = {2024},
doi = {10.1142/S0219199723500219},

URL = { 
    
        https://doi.org/10.1142/S0219199723500219
    
    

},
eprint = { 
    
        https://doi.org/10.1142/S0219199723500219
    
    

}
,}

@article{Ryan23,
author = {Hynd, Ryan},
title = {Evolution of Mixed Strategies in Monotone Games},
journal = {SIAM Journal on Optimization},
volume = {33},
number = {4},
pages = {2750-2771},
year = {2023},
doi = {10.1137/22M1486066},

URL = { 
    
        https://doi.org/10.1137/22M1486066
    
    

},
eprint = { 
    
        https://doi.org/10.1137/22M1486066
    
    

}
}

@article{pass2023generalized,
  title={Generalized barycenters and variance maximization on metric spaces},
  author={Pass, Brendan},
  journal={Journal of Fixed Point Theory and Applications},
  volume={25},
  number={1},
  pages={5},
  year={2023},
issn={1661-7746},
doi={10.1007/s11784-022-01015-x},
url={https://doi.org/10.1007/s11784-022-01015-x}
}

@article{BAJGIRAN2022111608,
title = {Uncertainty quantification of the 4th kind; optimal posterior accuracy-uncertainty tradeoff with the minimum enclosing ball},
journal = {Journal of Computational Physics},
volume = {471},
pages = {111608},
year = {2022},
issn = {0021-9991},
doi = {https://doi.org/10.1016/j.jcp.2022.111608},
url = {https://www.sciencedirect.com/science/article/pii/S0021999122006714},
author = {Hamed Hamze Bajgiran and Pau Batlle and Houman Owhadi and Mostafa Samir and Clint Scovel and Mahdy Shirdel and Michael Stanley and Peyman Tavallali},
}

@book{Brezis2011,
author="Br\'ezis, Haim",
title="Functional Analysis, Sobolev Spaces and Partial Differential Equations",
year="2011",
publisher="Springer New York",
address="New York, NY",
isbn="978-0-387-70914-7",
doi="10.1007/978-0-387-70914-7_3",
url="https://doi.org/10.1007/978-0-387-70914-7_3"
}

@article {MR287357,
    AUTHOR = {Crandall, M. G. and Liggett, T. M.},
     TITLE = {Generation of semi-groups of nonlinear transformations on
              general {B}anach spaces},
   JOURNAL = {Amer. J. Math.},
  FJOURNAL = {American Journal of Mathematics},
    VOLUME = {93},
      YEAR = {1971},
     PAGES = {265--298},
      ISSN = {0002-9327,1080-6377},
   MRCLASS = {47.50},
  MRNUMBER = {287357},
MRREVIEWER = {A.\ Pazy},
       DOI = {10.2307/2373376},
       URL = {https://doi.org/10.2307/2373376},
}

@article {Lu2022osrresolutionodeframeworkunderstanding,
    AUTHOR = {Lu, Haihao},
     TITLE = {An {$O(s^r)$}-resolution {ODE} framework for understanding
              discrete-time algorithms and applications to the linear
              convergence of minimax problems},
   JOURNAL = {Math. Program.},
  FJOURNAL = {Mathematical Programming},
    VOLUME = {194},
      YEAR = {2022},
    NUMBER = {1-2},
     PAGES = {1061--1112},
      ISSN = {0025-5610,1436-4646},
   MRCLASS = {90C26 (90C47)},
  MRNUMBER = {4445476},
MRREVIEWER = {Ke\ Ke\ Li},
       DOI = {10.1007/s10107-021-01669-4},
       URL = {https://doi.org/10.1007/s10107-021-01669-4},
}

@misc{cai2024convergenceminmaxlangevindynamics,
      title={Convergence of the Min-Max Langevin Dynamics and Algorithm for Zero-Sum Games}, 
      author={Yang Cai and Siddharth Mitra and Xiuyuan Wang and Andre Wibisono},
      year={2024},
      eprint={2412.20471},
      archivePrefix={arXiv},
      primaryClass={cs.GT},
      url={https://arxiv.org/abs/2412.20471}, 
}

@article{JKO98,
author = {Jordan, Richard and Kinderlehrer, David and Otto, Felix},
title = {The Variational Formulation of the Fokker--Planck Equation},
journal = {SIAM Journal on Mathematical Analysis},
volume = {29},
number = {1},
pages = {1--17},
year = {1998},
doi = {10.1137/S0036141096303359},

URL = { 
    
        https://doi.org/10.1137/S0036141096303359
    
    

},
eprint = { 
    
        https://doi.org/10.1137/S0036141096303359
    
    

}
,
}

@book{ArrowHurwiczUzawa58,
    AUTHOR = {Arrow, Kenneth J. and Hurwicz, Leonid and Uzawa, Hirofumi},
     TITLE = {Studies in linear and non-linear programming},
    SERIES = {Stanford Mathematical Studies in the Social Sciences},
    VOLUME = {II},
      NOTE = {With contributions by H. B. Chenery, S. M. Johnson, S. Karlin,
              T. Marschak, R. M. Solow},
 PUBLISHER = {Stanford University Press, Stanford, CA},
      YEAR = {1958},
     PAGES = {vii+229},
   MRCLASS = {90.00},
  MRNUMBER = {108399},
MRREVIEWER = {A.\ Charnes},
}

@book {Brezis73,
    AUTHOR = {Br\'ezis, H.},
     TITLE = {Op\'erateurs maximaux monotones et semi-groupes de
              contractions dans les espaces de {H}ilbert},
    SERIES = {North-Holland Mathematics Studies},
    VOLUME = {No. 5},
      NOTE = {Notas de Matem\'atica, No. 50. [Mathematical Notes]},
 PUBLISHER = {North-Holland Publishing Co., Amsterdam-London; American
              Elsevier Publishing Co., Inc., New York},
      YEAR = {1973},
     PAGES = {vi+183},
   MRCLASS = {47H05},
  MRNUMBER = {348562},
MRREVIEWER = {Bruce\ Calvert},
}

@article{STURM2005149,
title = {Convex functionals of probability measures and nonlinear diffusions on manifolds},
journal = {Journal de Mathématiques Pures et Appliquées},
volume = {84},
number = {2},
pages = {149--168},
year = {2005},
issn = {0021-7824},
doi = {https://doi.org/10.1016/j.matpur.2004.11.002},
url = {https://www.sciencedirect.com/science/article/pii/S0021782404001485},
author = {Karl-Theodor Sturm},
}

@article {MR772092,
    AUTHOR = {Bakry, Dominique and \'Emery, Michel},
     TITLE = {Hypercontractivit\'e{} de semi-groupes de diffusion},
   JOURNAL = {C. R. Acad. Sci. Paris S\'er. I Math.},
  FJOURNAL = {Comptes Rendus des S\'eances de l'Acad\'emie des Sciences.
              S\'erie I. Math\'ematique},
    VOLUME = {299},
      YEAR = {1984},
    NUMBER = {15},
     PAGES = {775--778},
      ISSN = {0249-6291},
   MRCLASS = {60J60 (58G32 60J35)},
  MRNUMBER = {772092},
MRREVIEWER = {Jacques\ Vauthier},
}

@article {NikaidoIsoda55,
    AUTHOR = {Nikaid\^o, Hukukane and Isoda, Kazuo},
     TITLE = {Note on non-cooperative convex games},
   JOURNAL = {Pacific J. Math.},
  FJOURNAL = {Pacific Journal of Mathematics},
    VOLUME = {5},
      YEAR = {1955},
     PAGES = {807--815},
      ISSN = {0030-8730,1945-5844},
   MRCLASS = {90.0X},
  MRNUMBER = {73910},
MRREVIEWER = {D.\ Gale},
       URL = {http://projecteuclid.org/euclid.pjm/1171984836},
}

@article{Benzi_Golub_Liesen_2005, title={Numerical solution of saddle point problems}, volume={14}, DOI={10.1017/S0962492904000212}, journal={Acta Numerica}, author={Benzi, Michele and Golub, Gene H. and Liesen, Jörg}, year={2005}, pages={1–137}}

@InProceedings{pmlr-v97-hsieh19b,
  title = 	 {Finding Mixed {N}ash Equilibria of Generative Adversarial Networks},
  author =       {Hsieh, Ya-Ping and Liu, Chen and Cevher, Volkan},
  booktitle = 	 {Proceedings of the 36th International Conference on Machine Learning},
  pages = 	 {2810--2819},
  year = 	 {2019},
  volume = 	 {97},
  url = 	 {https://proceedings.mlr.press/v97/hsieh19b.html}
}

@book{Bacak+2014,
url = {https://doi.org/10.1515/9783110361629},
title = {Convex Analysis and Optimization in Hadamard Spaces},
author = {Miroslav Bacak},
publisher = {De Gruyter},
address = {Berlin, München, Boston},
doi = {doi:10.1515/9783110361629},
isbn = {9783110361629},
year = {2014},
lastchecked = {2025-06-18}
}

@article{DeGiorgi1980,
  author    = {De Giorgi, Ennio and Marino, Antonio and Tosques, Mario},
  title     = {Problemi di evoluzione in spazi metrici e curve di massima pendenza},
  journal   = {Atti della Accademia Nazionale dei Lincei. Classe di Scienze Fisiche, Matematiche e Naturali. Rendiconti},
  volume    = {68},
  number    = {3},
  pages     = {180--187},
  month     = mar,
  year      = {1980},
  issn      = {0392-7881},
  publisher = {Accademia Nazionale dei Lincei},
  url       = {http://eudml.org/doc/288859},
  language  = {italian},
}
\end{document}